\documentclass[review,3p,square,sort]{elsarticle}

\usepackage{amssymb}
\usepackage{amsmath}
\usepackage{amsthm}
\usepackage{stmaryrd}
\usepackage{bbm}
\usepackage{enumitem}

\usepackage{fancyhdr}
\usepackage{hyperref}
\hypersetup{%
  pdftitle={Existence, uniqueness and stability of $L^1$ solutions for multidimensional BSDEs with generators of one-sided Osgood type},%
  pdfsubject={BSDEs with one-sided Osgood condition},%
  pdfauthor={ShengJun FAN},%
  pdfkeywords={BSDEs, One-sided Osgood condition, Existence and uniqueness, $L^1$ solution},%
  pdfstartview=FitH,
  CJKbookmarks=true,%
  bookmarksnumbered=true,%
  bookmarksopen=true,%
  colorlinks=true, linkcolor=blue, urlcolor=blue, citecolor=blue, %
}

\begin{document}
\journal{ArXiv}

 \newcommand{\eps}{\varepsilon}
 \newcommand{\lam}{\lambda}
 \newcommand{\To}{\rightarrow}
 \newcommand{\as}{{\rm d}\mathbbm{P}\times{\rm d}t-a.e.}
 \newcommand{\ass}{{\rm d}\mathbbm{P}\times{\rm d}s-a.e.}
 \newcommand{\ps}{{\rm d}\mathbbm{P}-a.s.}
 \newcommand{\jf}{\int_t^T}
 \newcommand{\tim}{\times}

 \newcommand{\F}{\mathcal{F}}
 \newcommand{\E}{\mathbb{E}}
 \newcommand{\N}{\mathbb{N}}
 \newcommand{\s}{\mathcal{S}}
 \newcommand{\M}{{\rm M}}
 \newcommand{\T}{[0,T]}
 \newcommand{\LT}{L^1(\Omega, \F_T,\mathbbm{P};\R^k)}
 \newcommand{\Lp}{L^p(\Omega, \F_T,\mathbbm{P};\R^k)}

 \newcommand{\R}{{\mathbb R}}
 \newcommand{\Q}{{\mathbb Q}}
 \newcommand{\RE}{\forall}

\newcommand {\Lim}{\lim\limits_{n\rightarrow \infty}}
\newcommand {\Dis}{\displaystyle}

\begin{frontmatter}
\title {Existence, uniqueness and stability of $L^1$ solutions for multidimensional BSDEs with generators of one-sided Osgood type\tnoteref{fund}}
\tnotetext[fund]{Supported by the National Natural Science Foundation of China (No. 11371362) and the Fundamental Research Funds for the Central Universities (No. 2012QNA36).\vspace{0.2cm}}

\author{ShengJun FAN}
\ead{f$\_$s$\_$j@126.com}

\address{School of Mathematics, China University of Mining and Technology, Xuzhou 221116, PR China\vspace{-0.8cm}}

\begin{abstract}
We establish a general existence and uniqueness result of $L^1$ solution for a multidimensional backward stochastic differential equation (BSDE for short) with generator $g$ satisfying a one-sided Osgood condition as well as a general growth condition in $y$, and a Lipschitz condition together with a sublinear growth condition in $z$, which improves some existing results. In particular, we put forward and prove a stability theorem of the $L^1$ solutions for the first time. A new type of $L^1$ solution is also investigated. Some delicate techniques involved in the relationship between convergence in $L^1$ and in probability and dividing appropriately the time interval play crucial roles in our proofs.\vspace{0.1cm}
\end{abstract}

\begin{keyword}
Backward stochastic differential equation \sep $L^1$ solution \sep Existence and uniqueness\sep Stability theorem\sep One-sided Osgood condition\vspace{0.3cm}

\MSC[2010] 60H10
\end{keyword}
\end{frontmatter}\vspace{-0.4cm}


\section{Introduction}

We fix a nonnegative real number $T$ as well as two positive integers $k$ and $d$, and let $\R^+:=[0,+\infty)$. Let $\mathbbm{1}_{A}$ represent the indicator function of a set $A$, and $\langle x,y\rangle$ the inner product of $x,y\in\R^k$. The Euclidean norms of a vector $y\in\R^k$ and a matrix $z\in \R^{k\times d}$ are defined by $|y|$ and $|z|$, respectively.

Assume that $(\Omega,\F,\mathbbm{P})$ is a completed probability space carrying a standard $d$-dimensional Brownian motion $(B_t)_{t\geq 0}$, and that $(\F_t)_{t\geq 0}$ is the natural $\sigma$-algebra filtration generated by $(B_t)_{t\geq 0}$ and $\F=\F_T$. For each $p>0$, denote by $\Lp$ the set of all $\R^k$-valued and $\F_T$-measurable random vectors $\xi$ such that $\E[|\xi|^p]<+\infty$, by ${\s}^p(0,T;\R^k)$ (or $\s^p$ simply) the set of $\R^k$-valued, $(\F_t)$-adapted and continuous processes $(Y_t)_{t\in\T}$ such that
$$\|Y\|_{{\s}^p}:=\left( \E\left[\sup_{t\in\T} |Y_t|^p\right] \right)^{1\wedge 1/p}<+\infty,$$
and by ${\rm M}^p(0,T;\R^{k\times d})$ (or ${\rm M}^p$ simply) the set of $(\F_t)$-progressively measurable ${\R}^{k\times d}$-valued processes $(Z_t)_{ t\in\T}$ such that
$$\|Z\|_{{\rm M}^p}:=\left\{ \E\left[\left(\int_0^T |Z_t|^2\ {\rm d}t\right)^{p/2}\right] \right\}^{1\wedge 1/p}<+\infty.
$$
It is well known that for each $p\geq 1$, $\s^p$ and ${\rm M}^p$ are both Banach spaces respectively endowed with the norms $\|\cdot\|_{{\s}^p}$ and $\|\cdot\|_{{\rm M}^p}$. And, for each $p\in (0,1)$, $\s^p$ and ${\rm M}^p$ are both complete metric spaces with the resulting distances $(Y,Y')\mapsto \|Y-Y'\|_{{\s}^p}$ and $(Z,Z')\mapsto \|Z-Z'\|_{{\rm M}^p}$ respectively.

We recall that a process $(Y_t)_{t\in\T}$ belongs to the class (D) if the family of variables $\{|Y_\tau|:\tau\ {\rm is\ an}$ $(\F_t){\rm - stopping\ time\ bounded\ by}\ n\}$
is uniformly integrable.

In this paper, we are interested in solving the following multidimensional backward stochastic differential equation (BSDE for short):
\begin{equation}
y_t=\xi+\int_t^Tg(s,y_s,z_s){\rm d}s-\int_t^Tz_s {\rm d}B_s,\ \   t\in\T,
\end{equation}
where $\xi\in\LT$ is called the terminal condition, $T$ is called the time horizon, and the random function $$g(\omega,t,y,z):\Omega\tim \T\tim {\R}^{k }\tim
{\R}^{k\times d}\longmapsto {\R}^k$$
is $(\F_t)$-progressively measurable for each $(y,z)$, called the generator of BSDE (1). Furthermore, the triple $(\xi,T,g)$ is usually called the parameters of BSDE (1).\vspace{0.1cm}

Throughout this paper, we use the following definitions on solutions of (1).\vspace{0.1cm}

{\bf Definition 1}\ \ A solution of BSDE (1) is a pair of $(\F_t)$-progressively measurable processes $(y_t,z_t)_{t\in\T}$ with values in ${\R}^k\times {\R}^{k\times d}$ such that $\ps$, $\int_0^T|z_t|^2\ {\rm d}t<+\infty$, $\int_0^T|g(t,y_t,z_t)|\ {\rm d}t<+\infty$, and (1) holds true for each $t\in\T$.\vspace{0.2cm}

{\bf Definition 2}\ \ Assume that $(y_t,z_t)_{t\in\T}$ is a solution of (1). If $(y_t,z_t)_{t\in\T}\in {\s}^p(0,T;\R^{k})\times {\rm M}^p(0,T;\R^{k\times d})$ for some $p>1$, then it is called an $L^p$ solution of BSDE (1); if $(y_t)_{t\in\T}$ belongs to the class (D) and
$(y_t,z_t)_{t\in\T}\in {\s}^\beta(0,T;\R^{k})\times {\rm M}^\beta(0,T;\R^{k\times d})$ for each $\beta\in (0,1)$, then it is called an $L^1$ solution of BSDE (1).

It is well known that nonlinear BSDEs were initially introduced in 1990 by \citet{Par90}. They put forward and proved an existence and uniqueness result for $L^2$ solution of multidimensional BSDEs under the Lipschitz assumption of $g$ as well as the square integrability assumption of $\xi$ and $g(t,0,0)$. From then on, the BSDE theory has attracted more and more interests, and due to the closely connections with many questions, it has gradually become a very powerful tool in many fields including stochastic control, financial mathematics, nonlinear mathematical expectation and partial differential equations, see \cite{Bah10,Buc00,Chen00,Delb10,El97,Hu11,Hu08,Jia10,Kob00,
Mor09,Par99,Peng97,Tang98,Xing12} and so on.

There is no doubt that the existence and uniqueness of the solution is one of the most fundamental and kernel problems in the study on the theory and application of BSDEs. From the beginning, many researchers have attempted to improve the result of the $L^2$ solution of \cite{Par90} by weakening the Lipschitz hypothesis on $g$, see, for example, \cite{Bri07,Chen00,El97,Fan10,
Fan13,FJD13,Ham03,Jia08,Jia10,Lep97,Mao95,Xu15} for a survey.
At the same time, the existence and uniqueness of the $L^p\ (p>1)$ solution for BSDEs has been extensively investigated by \cite{Bri03,El97,Fan15,FJ14}, etc. Starting around 1998, the existence and uniqueness of the bounded solution and the solution whose exponential moments of certain order exist have also been becoming one of emphasis in the study on BSDE theory, one can see \cite{Bri06,Bri08,Delb11,Hu15,Kob00,Lep98,Mor09,Rich12} for this topic, where the generator $g$ may have a quadratic or superquadratic growth in $z$.

On the other hand, in 1997, \citet{Peng97} introduced the notion of $g$-martingales by solutions of BSDEs, which can be viewed as some kind of nonlinear martingales. Since the classical theory of martingales is carried in the integrable space, the question of solving a BSDE with only integrable parameters comes up naturally. In this spirit, some recent works including \cite{Bri02,Bri03,Bri06,Fan12,FL10,Peng97,Xiao12,Xing12} investigated the existence and uniqueness of the $L^1$ solution of BSDEs. In particular, we would like to mention that \citet{Bri03} established a general existence and uniqueness result of $L^1$ solution for a multidimensional BSDE with generator $g$ satisfying a monotonicity condition (see (H1) with $\rho(x)=\mu x$ in Section 3) as well as a general growth condition in $y$ (see (H2) in Section 3), and a Lipschitz condition together with a sublinear growth condition in $z$ (see (H3) in Section 3). Here, we also mention that multidimensional BSDEs are more difficult to handle than the one-dimensional case since for multidimensional BSDEs we usually can not establish or employ the comparison theorem of solutions. And, it is well known that the $L^1$ solution is more difficult to treat than the $L^p\ (p>1)$ solution.

The present paper focus on the $L^1$ solution of multidimensional BSDEs. First of all, we will extend the existence and uniqueness result of the $L^1$ solution established in \cite{Bri03} by weakening the monotonicity condition of $g$ in $y$ to a one-sided Osgood condition (see (H1) in Section 3). Under a Osgood condition of $g$ in $y$ and a Lipschitz condition of $g$ in $z$, \citet{FJD13} first proved the existence and uniqueness of $L^2$ solution for multidimensional BSDEs. Recently, \citet{Fan15} further extended this result and established the existence and uniqueness of $L^p\ (p>1)$ solution for a multidimensional BSDE with generator $g$ satisfying a $p$-order weak monotonicity condition (see (H1a)$_p$ in Section 2) as well as a general growth condition in $y$, and a Lipschitz condition in $z$. We point out that in the case of $p=1$, the $p$-order weak monotonicity condition used in \cite{Fan15} becomes the one-sided Osgood condition used in this paper. From this point of view, it is very natural to investigate the $L^1$ solution of multidimensional BSDEs under the one-sided Osgood condition of $g$ in $y$. However, when we use this condition instead of the usual monotonicity condition, some essential difficulty arise especially in the proof of existence of the $L^1$ solution. By virtue of Gronwall's inequality, Bihari's inequality and the relationship between convergence in $L^1$ and in probability together with two updated apriori estimates established in \citet{Fan15}, we first consider the case when the generator $g$ is independent of $z$ (see Proposition 5 in Section 3). Then, making use of two estimates established in \citet{Xu15} and \citet{FJ14} together with Bihari's inequality, by a delicate argument involved in a Piciard's iterative procedure and a technique dividing the time interval $\T$ we prove the existence of the $L^1$ solution for the general case (see Theorem 2 in Section 3). Here, we mention that it is interesting that in the case of the $\alpha$ defined in (H3) values in $[1/2,1)$, the one-sided Osgood condition need to be replaced with a $p$-order ($p>1$) one-sided Mao's condition (see (H1b)$_p$ in Section 2).

The second objective of this paper is to put forward and prove a stability theorem for the $L^1$ solutions of multidimensional BSDEs with generators of one-sided Osgood type. To the best of our knowledge, this is the first time for the $L^1$ solution of multidimensional BSDEs. It is not very hard to obtain a stability result of $L^p\ (p>1)$ solutions for multidimensional BSDEs since by classical techniques one can establish and employ apriori estimates on the $L^p$ solution when $p>1$ (see, for example, Theorem 2 in \citet{Fan15} for more details). However, it is well known that when $p=1$ the apriori estimates with respect to the first component of the $L^1$ solution are not valid any longer especially when the generator $g$ depends on $z$, which brings intrinsic difficulty when one tries to establish the stability of $L^1$ solutions. This may be the reason that by far there is still no reported work on the stability of $L^1$ solutions for multidimensional BSDEs even when $g$ only satisfies the monotonicity condition or the uniformly Lipschitz condition in $y$ other than the one-sided Osgood condition. In this paper, we will fill up the gap. More specifically, enlightened by the proof of the existence of the $L^1$ solution in this paper, we will first introduce some auxiliary BSDEs by virtue of a Picard's iterative procedure as a bridge and then use a very delicate argument to establish a stability theorem of the $L^1$ solutions for multidimensional BSDEs with generators of one-sided Osgood type (see Theorem 4 in Section 5), where expect for Gronwall's inequality, Bihari's inequality, the relationship between convergence in $L^1$ and in probability, and the technique dividing the time interval, the induction technique and the sharp apriori estimates established in \cite{Fan15}, \cite{Xu15} and \cite{FJ14} all play important roles.

In addition, in this paper we also investigate, for the first time, the existence and uniqueness together with the stability of the solutions in the space $\s^1\times {\rm M}^1$ for multidimensional BSDEs (see Theorem 3 in Section 4 and Theorem 5 in Section 5).

The remainder of this paper is organized as follows. In Section 2 we gather several updated apriori estimates with respect to the solutions of multidimensional BSDEs and two technical lemmas. In section 3 we state and prove the existence and uniqueness result of $L^1$ solutions for the multidimensional BSDEs, and in Section 4 we are interested in solving the multidimensional BSDEs in $\s^1\times {\rm M}^1$ and provide two examples to illustrate our theoretical results. Finally, in Section 5 we put forward and prove a stability theorem of $L^1$ solutions as well as solutions in $\s^1\times {\rm M}^1$ for the multidimensional BSDEs.

\section{Preliminaries}

In this section, we first introduce several sharp apriori estimates with respect to solutions of multidimensional BSDEs, which will play very important roles in the proof of our main results. For this, let us introduce the following assumptions with respect to the generator $g$:\vspace{0.3cm}

{\bf (A1)}\ \ $\as,\ \RE\ (y,z)\in \R^k\times\R^{k\times d},\ \ \left\langle y,
g(\omega,t,y,z)\right\rangle\leq \mu |y|^2+\nu |y||z|+|y|f_t+\varphi_t,\ \ $\vspace{0.3cm}\\
where $\mu$ and $\nu$ are two positive constants, $f_t$ and $\varphi_t$ are two $(\F_t)$-progressively measurable and nonnegative processes satisfying
$$
\E\left[\left(\int_0^T f_t\ {\rm d}t\right)^p\right]<+\infty\ \ \ {\rm and}\ \ \ \E\left[\left(\int_0^T \varphi_t\ {\rm d}t\right)^{p/2}\right]<+\infty.\vspace{0.2cm}
$$

{\bf (A2)}\ \ $\as,\ \RE\ (y,z)\in \R^k\times\R^{k\times d}$,\vspace{0.25cm}\\
\hspace*{3cm}$|y|^{p-1}\left\langle {y\over |y|}\mathbbm{1}_{|y|\neq 0},
g(\omega,t,y,z)\right\rangle\leq \psi(|y|^p)+\nu |y|^{p-1}|z|+|y|^{p-1}f_t,$\vspace{0.35cm}\\
where $\nu>0$ is a constant, $f_t$ is an $(\F_t)$-progressively measurable and nonnegative process satisfying
$$\E\left[\left(\int_0^T f_t\ {\rm d}t\right)^p\right]<+\infty,$$
and $\psi(\cdot):\R^+\mapsto \R^+$ is a nondecreasing and concave function with $\psi(0)=0$.\vspace{0.3cm}

{\bf (A3)}\ \ $\as,\RE\ (y,z)\in \R^k\times\R^{k\times d},\ \ \left\langle {y\over |y|}\mathbbm{1}_{|y|\neq 0}, g(\omega,t,y,z)\right\rangle\leq
\phi^{{1\over p}}(|y|^p)+\nu|z|+f_t,\ \ $\vspace{0.3cm}\\
where $\nu>0$ is a constant, $f_t$ is an $(\F_t)$-progressively measurable and nonnegative process satisfying
$$\E\left[\left(\int_0^T f_t\ {\rm d}t\right)^p\right]<+\infty,$$
and $\phi(\cdot):\R^+\mapsto \R^+$ is a nondecreasing and concave function with $\phi(0)=0$.\vspace{0.2cm}

The above assumptions (A2) and (A3) are respectively related to  following assumptions (H1a)$_p$ and (H1b)$_p$, which are put forward and used in \citet{Fan15} at the first time. Assumptions (H1a)$_p$ and (H1b)$_p$ will also be employed in this paper.\vspace{0.2cm}

{\bf (H1a)$_p$} $g$ satisfies a $p$-order weak monotonicity condition in $y$, i.e., there exists a nondecreasing and concave function $\kappa(\cdot):\R^+\mapsto \R^+$ with $\kappa(0)=0$, $\kappa(u)>0$ for $u>0$ and $\int_{0^+} {{\rm d}u\over \kappa(u)}=+\infty$ such that $\as$, $\RE\ y_1,y_2\in \R^k,z\in\R^{k\times d}$,
$$
|y_1-y_2|^{p-1}\langle {y_1-y_2\over |y_1-y_2|}\mathbbm{1}_{|y_1-y_2|\neq 0},g(\omega,t,y_1,z)-g(\omega,t,y_2,z)\rangle\leq \kappa(|y_1-y_2|^p),
$$
where and hereafter, $\int_{0^+} {{\rm d}u\over \kappa(u)}:=\lim\limits_{\epsilon \rightarrow 0}\int_{0}^{\epsilon} {{\rm d}u\over \kappa(u)};\vspace{0.2cm}$

{\bf (H1b)$_p$} $g$ satisfies a $p$-order one-sided Mao's condition in $y$, i.e., there exists a nondecreasing and concave function $\varrho(\cdot):\R^+\mapsto \R^+$ with $\varrho(0)=0$, $\varrho(u)>0$ for $u>0$ and $\int_{0^+} {{\rm d}u\over \varrho(u)}=+\infty$ such that $\as$, $\RE\ y_1,y_2\in \R^k,z\in\R^{k\times d}$,
$$
\langle {y_1-y_2\over |y_1-y_2|}\mathbbm{1}_{|y_1-y_2|\neq 0},
g(\omega,t,y_1,z)-g(\omega,t,y_2,z)\rangle\leq \varrho^{{1\over p}}(|y_1-y_2|^p).\vspace{0.2cm}
$$

The following Propositions 1-2 are respectively Propositions 2-3 in \citet{Fan15}, and the following Proposition 3 comes from Proposition 1 in \citet{Xu15}.\vspace{0.1cm}

{\bf Proposition 1} Let $p>0$ and (A1) hold. Suppose that $(y_t,z_t)_{t\in\T}$ is a solution of BSDE (1) such that $y_\cdot\in {\s}^p(0,T;\R^{k})$. Then $z_\cdot$ belongs to ${\rm M}^p(0,T;\R^{k\times d})$, and $\ps$, for each $0\leq u\leq t\leq T$,
$$
\begin{array}{lll}
\Dis \E\left[\left.\left(\int_t^T |z_s|^2\ {\rm
d}s\right)^{p/2}\right|\F_u\right]&\leq & \Dis
C_{\mu,\nu,p,T}\E\left[\left.\sup\limits_{s\in
[t,T]}|y_s|^p\right|\F_u\right]+C_p\E\left[\left.
\left(\int_t^T f_s\ {\rm d}s\right)^p\right|\F_u\right]\vspace{0.1cm}\\
&& \Dis+C_p\E\left[\left.\left(\int_t^T \varphi_s\ {\rm d}s\right)^{p/2}\right|\F_u\right],
\end{array}
$$
where $C_{\mu,\nu,p,T}>0$ is a constant depending on $(\mu,\nu,p,T)$, and $C_p>0$ is a constant depending only on $p$.\vspace{0.3cm}

{\bf Proposition 2} Let $p>1$ and (A2) hold. Suppose that $(y_t,z_t)_{t\in\T}$ is an $L^p$ solution of BSDE (1). Then, there exists a constant $C_{\nu,p}>0$ depending only on $\nu, p$ such that $\ps$, for each $0\leq u\leq t\leq T$,\vspace{0.2cm}
$$
\E\left[\left.\sup\limits_{s\in
[t,T]}|y_s|^{p}\right|\F_u\right]\leq
e^{C_{\nu,p}(T-t)}\left\{\E[\left.|\xi|^p\right|\F_u]+\int_t^T
\psi(\E[\left.|y_s|^p\right|\F_u]) \ {\rm
d}s+\E\left[\left.\left(\int_t^T f_s\ {\rm
d}s\right)^{p}\right|\F_u\right]\right\}.\vspace{0.3cm}
$$

{\bf Proposition 3} Let $g$ satisfy (A2) with $p=2$.  Suppose that $(y_t,z_t)_{t\in\T}$ is an $L^2$ solution of BSDE (1). Then, there exists a constant $C_{\nu}>0$ depending only on $\nu$ such that $\ps$, for each $0\leq u\leq t\leq T$,\vspace{0.1cm}
$$
\begin{array}{lll}
&&\Dis \E\left[\left.\sup\limits_{r\in
[t,T]}|y_r|^2\right|\F_u\right]+\E\left[\left.\int_t^T
|z_s|^2\ {\rm d}s\right|\F_u\right]\vspace{0.1cm}\\
&\leq & \Dis
e^{C_{\nu}(T-t)}\left\{\E\left[\left.|\xi|^2\right|\F_u\right]+
\int_t^T
\psi\left(\E\left[\left.|y_s|^2\right|\F_u\right]\right)\ {\rm d}s+\E\left[\left.\left(\int_t^T f_s\ {\rm
d}s\right)^2\right|\F_u\right]\right\}.
\end{array}\vspace{0.2cm}
$$

In the same way as that in Lemmas 2-3 of \citet{FJ14}, we can prove the following proposition. The proof is omitted here.\vspace{0.1cm}

{\bf Propositions 4}\ Let $p>1$ and (A3) hold. Suppose that  $(y_t,z_t)_{t\in\T}$ is a solution of BSDE (1) such that $y_\cdot\in {\s}^p(0,T;\R^{k})$. Then, $z_\cdot$ belongs to ${\rm M}^p(0,T;\R^{k\times d})$ and there exists a positive constant $C_{\nu,p}$ depending on $\nu$ and $p$ such that $\ps$, for each $0\leq u\leq t\leq T$,\vspace{0.1cm}
$$
\begin{array}{lll}
&&\Dis \E\left[\left.\sup\limits_{s\in
[t,T]}|y_s|^{p}\right|\F_u\right]+\E\left[\left.
\left(\int_t^T|z_s|^2\ {\rm d}s\right)^{p/2}\right|\F_u\right]\vspace{0.1cm}\\
&\leq & \Dis e^{C_{\nu,p}(T-t)}\left\{\E\left[\left.|\xi|^p\right|
\F_u\right]+\int_t^T
\phi\left(\E\left[\left.|y_s|^p\right|\F_u\right]\right) \ {\rm d}s+\E\left[\left.\left(\int_t^T f_s\ {\rm
d}s\right)^p\right|\F_u\right]\right\}.
\end{array}\vspace{0.3cm}
$$

Now, let us introduce two technical lemmas, which will be used later.  Firstly, the following Lemma 1 gives a sequence of upper bounds for linear growth functions, which comes from \citet{Fan10}.\vspace{0.2cm}

{\bf Lemma 1}\ Suppose that $\bar\kappa(\cdot):\R^+\mapsto \R^+$ increases at most linearly, i.e., there exists a constant $A>0$ such that $$\bar\kappa(x)\leq A(x+1),\ \ \RE\ x\in \R^+.$$
Then for each $m\geq 1$, we have
$$
\bar\kappa(x)\leq (m+2A)x+\bar\kappa\left({2A\over m+2A}\right),\ \ \RE\ x\in \R^+.\vspace{0.2cm}
$$

The following Lemma 2 can be regarded as a backward version of classical Bihari's inequality, which can be proved by classical methods. The proof is omitted.\vspace{0.2cm}

{\bf Lemma 2}\ (Bihari's inequality)\ Let the nonnegative function $u(\cdot):[0,T]\mapsto \R^+$ satisfy
$$
u(t)\leq u_0 +\int_t^T \bar\psi(u(s))\ {\rm d}s,\ \ t\in \T,
$$
where $u_0$ is a positive real number, $\bar\psi(\cdot):\R^+\mapsto \R^+$ is a continuous and nondecreasing function, $\bar\psi(0)=0$, $\bar\psi(u)>0$ for $u>0$ and $\int_{0^+} {1\over \bar\psi(u)}\ {\rm d}u=+\infty$. Then, for each $t\in \T$, we have
$$
u(t)\leq \Psi^{-1}(\Psi(u_0)+T-t),
$$
where
$$
\Psi(x):=\int_1^x {1\over \bar\psi(u)}\ {\rm d}u, \ \ x>0\vspace{0.2cm}
$$
is a strictly increasing function valued in $\R$, and $\Psi^{-1}$ is the inverse function of $\Psi$. In particular, if $u_0=0$, then $u(t)=0$ for each $t\in \T$.\vspace{0.2cm}

To the end of this section, we would like to especially mention that even though Propositions 1-4 mentioned above appear similar, there are some distinguish differences among both their conditions and conclusions. They will also play different roles in the proof of our main results.

\section{Existence and uniqueness of the $L^1$ solution}

We first introduce the following assumptions
on the generator $g$ used in \citet{Fan15}, \citet{Fan13} and \citet{Bri03}:\vspace{0.2cm}

{\bf (H1)} $g$ satisfies a one-sided Osgood condition in $y$, i.e., there exists a nondecreasing and concave function $\rho(\cdot):\R^+\mapsto \R^+$ with $\rho(0)=0$, $\rho(u)>0$ for $u>0$ and $\int_{0^+} {{\rm d}u\over \rho(u)}=+\infty$ such that $\as$,
$$
\RE\ y_1,y_2\in \R^k,z\in\R^{k\times d},\ \ \left\langle {y_1-y_2\over |y_1-y_2|}\mathbbm{1}_{|y_1-y_2|\neq 0},
g(\omega,t,y_1,z)-g(\omega,t,y_2,z)\right\rangle\leq \rho(|y_1-y_2|).\vspace{0.2cm}
$$

{\bf (H2)} $g$ has a general growth with respect to $y$, i.e, \vspace{0.2cm}
$$\RE\ r>0,\ \ \E\left[\int_0^T \bar\phi_r(t)\ {\rm d}t\right]<+\infty\ \ {\rm with}\ \ \bar\phi_r(t):=\sup\limits_{|y|\leq r}
|g(\omega,t,y,0)|;$$
Furthermore, $\as$, $\RE\ z\in {\R^{k\times d}},\ \ \ y\longmapsto g(\omega,t,y,z)$ is continuous.\vspace{0.2cm}

{\bf (H3)}\ $g$ is Lipschitz continuous in $z$, uniformly with respect to $(\omega,t,y)$, i.e., there exists a constant $\lam\geq 0$ such that $\as$,
$$\RE\ y\in \R^k,z_1,z_2\in\R^{k\times d},\ \ |g(\omega,t,y,z_1)-g(\omega,t,y,z_2)|\leq \lam |z_1-z_2|;$$
Furthermore, $g$ has a sublinear growth in $z$, i.e., there exist two constants $\gamma>0$ and $\alpha\in (0,1)$ as well as an $(\F_t)$-progressively measurable and nonnegative process $(g_t)_{t\in\T}$ satisfying $\E\left[\int_0^Tg_t\ {\rm d}t\right]<+\infty$ such that $\as$, \vspace{0.05cm}
$$
\RE\ y\in \R^k,z\in\R^{k\times d},\ \ |g(\omega,t,y,z)-g(\omega,t,y,0)|\leq \gamma(g_t(\omega)+|y|+|z|)^\alpha.\vspace{0.2cm}
$$

{\bf Remark 1} For later use, it follows from Proposition 1 in \citet{Fan15} that for each $p>1$,
$${\rm (H1b)}_p \Longrightarrow {\rm (H1)}\Longrightarrow {\rm (H1a)}_p,$$
and when $p=1$, they are same. In addition, the functions $\rho(\cdot)$, $\kappa(\cdot)$ and $\varrho(\cdot)$ in (H1), (H1a) and (H1b) all increase at most linearly since they are all nondecreasing and concave function valued $0$ at $0$. Here and hereafter we will always denote by $A$ the linear-growth constants of them, i.e.,
$$\rho(x)\leq A(x+1),\ \ \kappa(x)\leq A(x+1),\ \ \varrho(x)\leq A(x+1),\ \ \RE\ x\in\R^+.$$
Finally, by Proposition 1 in \citet{Fan15} we also point out that the concavity condition of $\rho(\cdot)$ and $\varrho(\cdot)$ defined respectively in assumptions (H1) and (H1b)${}_p$ can be replaced with the continuity condition.\vspace{0.2cm}

The following Theorems 1-2 are the main results of this section.\vspace{0.1cm}

{\bf Theorem 1} Assume that the generator $g$ satisfies assumptions (H1) and (H3). Then for each $\xi\in\LT$, BSDE (1) admits at most one solution $(y_\cdot,z_\cdot)$ such that $y_\cdot$ belongs to the class (D) and $z_\cdot$ belongs to $\bigcup_{\beta>\alpha}\M^\beta$,
which leads to that it admits at most one $L^1$ solution.\vspace{0.1cm}

{\bf Proof}\ \ Assume that (H1) and (H3) hold and that both $(y_t,z_t)_{t\in\T}$ and $(y'_t,z'_t)_{t\in\T}$ are solutions of BSDE (1) such that both $(y_t)_{t\in\T}$ and $(y'_t)_{t\in\T}$ belong to the class (D), and both $(z_t)_{t\in\T}$ and $(z'_t)_{t\in\T}$ belong to $\M^\beta$ for some $\beta\in(\alpha, 1)$.\vspace{0.1cm}

We first show that $(y_t-y'_t)_{t\in\T}\in \s^{\beta/\alpha}$. In fact, let us fix $n\geq 1$ and denote $\tau_n$ the stopping time
$$
\tau_n:=\inf\left\{t\in\T:\int^t_{0}(|z_s|^2+|z'_s|^2)
\,\mathrm{d} s\geq n\right\}\wedge T.
$$
Corollary 2.3 in \citet{Bri03} leads to the following inequality with setting $\hat{y}_\cdot:=y_\cdot-y'_\cdot$ and   $\hat{z}_\cdot:=z_\cdot-z'_\cdot$, and $t\in\T$,
\begin{equation}
\Dis|\hat{y}_{t\wedge\tau_n}|\leq \Dis |\hat{y}_{\tau_n}|+\int^{\tau_n}_{t\wedge\tau_n}
\left\langle {\hat{y}_s\over |\hat{y}_s| }\mathbbm{1}_{|\hat{y}_s|\neq 0},g(s,y_s,z_s)-g(s,y'_s,z'_s)\right\rangle\
\mathrm{d}s -\int^{\tau_n}_{t\wedge\tau_n}\left\langle {\hat{y}_s\over |\hat{y}_s|}\mathbbm{1}_{|\hat{y}_s|\neq 0},\hat{z}_s\mathrm{d}B_s\right\rangle.
\end{equation}
It follows from assumptions (H1) and (H3) that $\ass$,
\begin{equation}
\begin{array}{ll}
&\Dis \left\langle{\hat{y}_s\over |\hat{y}_s|} \mathbbm{1}_{|\hat{y}_s|\neq 0},
g(s,y_s,z_s)-g(s,y'_s,z'_s)\right\rangle\vspace{0.1cm}\\
\leq &\Dis \left\langle{\hat{y}_s\over |\hat{y}_s|} \mathbbm{1}_{|\hat{y}_s|\neq 0},
g(s,y_s,z_s)-g(s,y'_s,z_s)\right\rangle
+|g(s,y'_s,z_s)-g(s,y'_s,z'_s)|\vspace{0.1cm}\\
\leq &\Dis\rho(|\hat{y}_s|)
+2\gamma\left(g_s+|y'_s|+|z'_s|+|z_s|\right)^\alpha.
\end{array}
\end{equation}
Then, combining (2) with (3) we can deduce that for each $n\geq 1$ and $t\in\T$,
\begin{equation}
|\hat{y}_{t\wedge\tau_n}|\leq\E\left[|\hat{y}_{\tau_n}|
+\left.\int^{\tau_n}_{t\wedge\tau_n}\rho(|\hat{y}_s|)\ \mathrm{d}s\right|\F_t\right]+G(t),
\end{equation}
where
$$
G(t):=2\gamma\E\left[\left.\int^{T}_{0}\left(g_s+
|y'_s|+|z'_s|+|z_s|\right)^\alpha\ \mathrm{d}s\right|\F_t\right].
$$
Furthermore, since $y'_\cdot$ belongs to the class (D), both $z_\cdot$ and $z'_\cdot$ belong to ${\rm M}^\beta$ with $\beta>\alpha$, and $\E\left[\int_0^Tg_t\ {\rm d}t\right]<+\infty$, we can use Doob's inequality, H\"{o}lder's inequality and Jensen's inequality to obtain that
\begin{equation}
\E\left[\sup_{t\in\T}|G(t)|^{\beta/\alpha}\right]<+\infty.
\end{equation}
Thus, since $\hat y_\cdot$ belongs to the class (D) and $\rho(\cdot)$ increases at most linearly, we can send $n$ to $+\infty$ in (4) and use Lebesgue's dominated convergence theorem, in view of $\tau_n\To T$ as $n\To \infty$, $\hat{y}_T=0$ and Remark 1, to get that for each $t\in \T$,
$$
|\hat{y}_t|\leq G(t)+\E\left[\left.\int_t^T
\rho(|\hat{y}_s|)\ \mathrm{d}s\right|\F_t\right]\leq AT+G(t)+A\int_t^T\E\left[\left.
|\hat{y}_s|\right|\F_t\right]\mathrm{d}s,
$$
and then
$$
\E\left[\left.|\hat{y}_r|\right|\F_t\right]\leq AT+G(t)+A\int_r^T
\E\left[\left.|\hat{y}_s|\right|\F_t\right]{\rm d}s,\quad r\in [t,T].\vspace{0.2cm}
$$
Gronwall's inequality yields that $\E\left[\left.|\hat{y}_r|\right|\F_t\right]\leq
\left(AT+G(t)\right)\cdot e^{A(T-r)},\ r\in [t,T]$, form which, by letting $r=t$, we have
\begin{equation}
|\hat{y}_t|\leq \left(AT+G(t)\right)\cdot e^{AT}.
\end{equation}
This inequality together with (5) leads to
\begin{equation}
\hat y_\cdot=(y_\cdot-y'_\cdot)\in \s^p\ \text{with}\ p:=\beta/\alpha>1.
\end{equation}

In the sequel, note that $(\hat y_t, \hat z_t)_{t\in\T}$ is a solution of the following BSDE:
\begin{equation}
\hat y_t=\int^T_t\hat g(s,\hat y_s,\hat z_s)\ \mathrm{d} s-\int^T_t\hat z_s\,\mathrm{d} B_s,\quad t\in\T,
\end{equation}
where for each $(y,z)\in\R^k\times\R^{k\times d}$,
$\hat g(t,y,z):=g(t,y+y'_t,z+z'_t)
-g(t,y'_t,z'_t)$.
It follows from assumptions (H1) and (H3) on $g$ together with Remark 1 that $\as$, for each $(y,z)\in \R^k\times\R^{k\times d}$,
\begin{equation}
\begin{array}{lll}
\left\langle y,\hat g(t,y,z)\right\rangle&\leq&\Dis \left\langle y, g(t,y+y'_t,z+z'_t)
-g(t,y'_t,z+z'_t)\right\rangle+|y||g(t,y'_t,z+z'_t)-g(t,y'_t,z'_t)|\\
&\leq&\Dis \bar\kappa(|y|^2)+\lam |y||z|\leq \Dis A|y|^2+\lam |y||z|+A,
\end{array}
\end{equation}
and
\begin{equation}
\begin{array}{lll}
\Dis |y|^{p-1}\left\langle {y\over |y|}\mathbbm{1}_{|y|\neq 0},\hat g(t,y,z)\right\rangle
&\leq&\Dis |y|^{p-1}\left\langle {y\over |y|}\mathbbm{1}_{|y|\neq 0}, g(t,y+y'_t,z+z'_t)
-g(t,y'_t,z+z'_t)\right\rangle\\
&&\Dis +|y|^{p-1}|g(t,y'_t,z+z'_t)-g(t,y'_t,z'_t)|
\vspace{0.1cm} \\
&\leq&\Dis \kappa(|y|^p)+\lam |y|^{p-1}|z|,
\end{array}
\end{equation}
where the functions $\bar\kappa(\cdot)$ and $\kappa(\cdot)$ are respectively defined in (H1a)$_2$ and (H1a)$_p$. Thus, on one hand, inequality (9) means that the generator $\hat g$ of BSDE (8) satisfies assumption (A1) with $\mu=A,\ \nu=\lam,\ f_t\equiv 0\ \text{and}\ \varphi_t\equiv A$. It then follows from Proposition 1 together with (7) that $(\hat y_t, \hat z_t)_{t\in\T}$ is an $L^p$ solution of BSDE (8). On the other hand, inequality (10) means that the generator $\hat g$ of BSDE (8) also satisfies assumption (A2) with $\psi(\cdot)=\kappa(\cdot),\ \nu=\lam\ \text{and}\ f_t\equiv 0$. It then follows from Proposition 2 with $u=0$ that there exists a positive constant $C_{\lam,p,T}$ depending only on $\lam$, $p$ and $T$ such that for each $t\in \T$,\vspace{0.1cm}
\begin{equation}
\E\left[\sup_{s\in [t,T]}|\hat y_s|^p\right]\leq C_{\lam,p,T} \int_t^T\kappa\left(\E\left[|\hat y_s|^p\right]\right)\ {\rm d}s\leq C_{\lam,p,T} \int_t^T\kappa\left(\E\left[\sup_{u\in [s,T]}|\hat y_u|^p\right]\right)\ {\rm d}s.\vspace{0.1cm}
\end{equation}
Thus, in view of the fact that $\int_{0^+} {{\rm d}u\over \kappa(u)}=+\infty$, Bihari's inequality (Lemma 2) yields that
\begin{equation}
\E\left[\sup_{t\in [0,T]}|y_t-y'_t|^p\right]=\E\left[\sup_{t\in [0,T]}|\hat y_t|^p\right]=0.\vspace{0.2cm}
\end{equation}

Finally, by (9), Remark 1 and Lemma 1 we can check that the generator $\hat g$ of BSDE (8) satisfies assumption (A1) with
$\mu=m+2A,\ \nu=\lam,\ f_t\equiv 0\ {\rm and}\ \ \varphi_t=\bar\kappa({2A\over m+2A})$
for each $m\geq 1$. It then follows from Proposition 1 with $u=t=0$ that there exists a positive constant $C_{m,\lam,p,T}$ depending on $m$, $\lam$, $p$ and $T$, and a positive constant $C_{p}$ depending only on $p$ such that for each $m\geq 1$,
\begin{equation}
\E\left[\left(\int_0^T |\hat z_s|^2\ {\rm
d}s\right)^{p/2}\right]\leq C_{m,\lam,p,T}\E\left[\sup\limits_{t\in [0,T]}|\hat y_t|^p\right]+C_{p}\left(\bar\kappa({2A\over m+2A})\cdot T\right)^{p/2}.
\end{equation}
Thus, in view of (12) and the fact that $\bar\kappa(\cdot)$ is a continuous function with $\bar\kappa(0)=0$, sending $m\To\infty$ in the previous inequality we deduce that
$$
\E\left[\left(\int_0^T |z_s-z'_s|^2\
{\rm d}s\right)^{p/2}\right]=\E\left[\left(\int_0^T |\hat z_s|^2\ {\rm
d}s\right)^{p/2}\right]=0.
$$
The proof of Theorem 1 is then complete. \vspace{0.3cm}\hfill $\Box$

We now turn to the existence part of our study. We will prove the following result.\vspace{0.1cm}

{\bf Theorem 2} Assume that the generator $g$ satisfies assumptions (H1)-(H3). In the case when the $\alpha$ defined in (H3) values in $[1/2,1)$, we further assume that there exists a constant $\bar p>1$ such that the function $\rho(\cdot)$ in (H1) satisfies
\begin{equation}
\int_{0^+}{u^{\bar p-1}\over \rho^{\bar p}(u)}{\rm d}u=+\infty.\vspace{0.2cm}
\end{equation}
Then for each $\xi\in\LT$, BSDE (1) admits an $L^1$ solution $(y_t,z_t)_{t\in\T}$.\vspace{0.3cm}

{\bf Remark 2} It follows from Proposition 1 in \citet{Fan15} that if $g$ satisfies (H1) with a $\rho(\cdot)$ satisfying (14) for some $\bar p>1$, then $g$ must satisfy (H1b)$_{\bar p}$. And vice versa. This fact will be perfectly utilized later.\vspace{0.2cm}

The following Proposition 5 is the first step to  prove Theorem 2, which studies the case where the generator $g$ does not depend on the variable $z$.\vspace{0.2cm}

{\bf Proposition 5} Let the generator $g$ be independent of $z$ and satisfy assumptions (H1)-(H2). Then for each $\xi\in\LT$, BSDE (1) admits an $L^1$ solution $(y_t,z_t)_{t\in\T}$.\vspace{0.2cm}

{\bf Proof}\ \ Assume that $\xi\in\LT$, $g$ is independent of $z$ and assumptions (H1) and (H2) hold. For each $n\geq 1$, we denote $q_n(x)=xn/(n\vee |x|)$ and set
$$
\xi^n:=q_n(\xi),\quad   g^n(t,y):=g(t,y)-g(t,0)+q_n(g(t,0)).
$$
Note that both $|\xi^n|$ and $|g^n(t,0)|$ are bounded by $n$, and that the generator $g^n(t,y)$ satisfies (H1) and  (H2) for each $n\geq 1$. It then follows from
Corollary 2 with $p=2$ in \citet{Fan15} that the following BSDE
\begin{equation}
y^n_t=\xi^n+\int^{T}_{t}g^n(s,y^n_s)\,\mathrm{d} s-\int^{T}_{t}z_s^n\,\mathrm{d} B_s,\quad t\in\T
\end{equation}
admits a unique $L^2$ solution $(y^n_t,z^n_t)_{t\in\T}$.\vspace{0.2cm}

For each $n,i\geq 1$, we set $\hat{y}^{n,i}_\cdot:=y^{n+i}_\cdot-y^n_\cdot$, $\hat z^{n,i}_\cdot:=z^{n+i}_\cdot-z^n_\cdot$ and $\hat\xi^{n,i}:=\xi^{n+i}-\xi^n$. It follows from Corollary 2.3 in \cite{Bri03} that for each $n,i\geq 1$ and $t\in\T$,
\begin{equation}
|\hat{y}^{n,i}_t|\leq \Dis |\hat{\xi}^{n,i}|+\int^{T}_{t}
\left\langle {\hat y^{n,i}_s\over |\hat y^{n,i}_s|}\mathbbm{1}_{|\hat y^{n,i}_s|\neq 0}, g^{n+i}(s,y^{n+i}_s)-g^n(s,y^n_s)\right\rangle
\ \mathrm{d} s\Dis -\int^{T}_{t}\left\langle {\hat y^{n,i}_s\over |\hat y^{n,i}_s|}\mathbbm{1}_{|\hat y^{n,i}_s|\neq 0},\hat z^{n,i}_s\,\mathrm{d} B_s\right\rangle.
\end{equation}
It follows from assumption (H1) and definition of $g^n(t,y)$ that $\ass$,
\begin{equation}
\begin{array}{lll}
\Dis \left\langle {\hat y^{n,i}_s\over |\hat y^{n,i}_s|}\mathbbm{1}_{|\hat y^{n,i}_s|\neq 0}, g^{n+i}(s,y^{n+i}_s)-g^n(s,y^n_s)\right\rangle
&\leq &\Dis \left\langle {\hat y^{n,i}_s\over |\hat y^{n,i}_s|}\mathbbm{1}_{|\hat y^{n,i}_s|\neq 0}, g^{n+i}(s,y^{n+i}_s)-g^{n+i}(s,y^n_s)\right\rangle
\vspace{0.1cm}\\
&&+|g^{n+i}(s,y^n_s)-g^n(s,y^n_s)|\vspace{0.1cm} \\
&\leq &\Dis \rho(|\hat y^{n,i}_s|)+|g(s,0)|       \mathbbm{1}_{|g(s,0)|>n}.
\end{array}
\end{equation}
Then, combining (16) and (17), in view of Fubini's Theorem and Jensen's inequality, we can get that for each $n,i\geq 1$ and $t\in\T$,
\begin{equation}
|\hat y_t^{n,i}|\leq H_n(t)+\int^{T}_{t}\rho\left(\E\left[\left.|\hat y^{n,i}_s|\right|\F_t\right]\right)\mathrm{d}s,
\end{equation}
where
$$H_n(t):=\E\left[|\xi|\mathbbm{1}_{|\xi|>n} +\left.\int^{T}_{0}|g(s,0)|\mathbbm{1}_{|g(s,0)|>n}
\mathrm{d}s\right|\F_t\right].\vspace{0.1cm}
$$

In the sequel, by virtue of Fubini's theorem and Jensen's inequality, it follows from (18) that for each $n,i\geq 1$ and $t\in\T$,
$$
\E\left[\left.|\hat y_r^{n,i}|\right|\F_t\right]\leq H_n(t)+\int^{T}_{r}\rho\left(\E\left[\left.|\hat y^{n,i}_s|\right|\F_t\right]\right)\mathrm{d}s,\quad r\in [t,T].
$$
Then in view of Remark 1, Gronwall's inequality yields that for each $n,i\geq 1$ and $(\omega,t)\in \Omega\times\T$,
\begin{equation}
|\hat{y}_t^{n,i}|=\E\left[\left.|\hat{y}_{t}^{n,i}|
\right|\F_t\right]\leq \left(H_n(t)+AT\right)\cdot e^{A(T-t)}.
\end{equation}
This inequality together with Lemma 6.1 in \cite{Bri03} leads to that for each $\beta\in (0,1)$,
\begin{equation}
\begin{array}{lll}
\Dis\sup\limits_{n\geq 1}\E\left[\sup\limits_{i\geq 1}\sup\limits_{t\in\T} |\hat{y}_t^{n,i}|^\beta\right]&\leq &\Dis C\left(1+\sup\limits_{n\geq 1} \E\left[\sup\limits_{t\in\T} |H_n(t)|^\beta\right]\right)\vspace{0.1cm}\\
&\leq & \Dis C\left(1+{1\over 1-\beta}\sup\limits_{n\geq 1}(\E[H_n(T)])^\beta\right)<+\infty,
\end{array}
\end{equation}
where $C>0$ is a constant depending only on $A,T$ and $\beta$. This means that the sequence
$$\left\{\sup_{i\geq 1}\sup_{t\in\T} |\hat{y}_t^{n,i}|^{\beta'}\right\}_{n=1}^{+\infty}
$$
is uniformly integrable for each $\beta'\in (0,1)$.\vspace{0.1cm}

On the other hand, by virtue of Fubini's theorem and Jensen's inequality, it follows from (18) that for each $n,i\geq 1$ and $t\in\T$,
$$
\E\left[\left.|\hat{y}_{t}^{n,i}|\right|\F_u\right]\leq
H_n(u)+\int_t^T\rho\left(\E\left[\left.|\hat{y}_{s}^{n,i}|
\right|\F_u \right]\right){\rm d}s,\ \ u\in[0,t],
$$
and then, in view of the fact that $\rho(\cdot)$ is a nondecreasing function,
\begin{equation}
\Dis h_n(t)\leq \sup_{u\in\T}H_n(u)+\int_t^T\rho\left(h_n(s)\right){\rm d}s,\ \ t\in\T,
\end{equation}
where
$$
h_n(t):=\sup_{i\geq 1}\sup\limits_{0\leq u\leq t}
\E\left[\left.|\hat{y}_{t}^{n,i}|
\right|\F_u\right]\vspace{0.2cm}
$$
are all nonnegative functions. Thus, by virtue of  Lemma 2 we can deduce that for each $n\geq 1$ and $\omega\in \Omega$,
\begin{equation}
h_n(t)\leq \Theta^{-1}(\Theta(\sup_{t\in\T}H_n(t))+T),\ \ t\in\T,
\end{equation}
where
\begin{equation}
\Theta(x):=\int_1^x {1\over \rho(x)}{\rm d}x,\ x>0 \vspace{0.2cm}
\end{equation}
is a strictly increasing and continuous function valued in $\R$, and $\Theta^{-1}$ is the inverse function of $\Theta$. Furthermore, by the maximum inequality with respect to sub-martingale and Lebesgue's dominated convergence theorem we have, for each $\epsilon>0$, as $n\To \infty$,
\begin{equation}
\mathbbm{P}\left(\left\{\sup_{t\in\T} |H_n(t)|\geq \epsilon\right\}\right)\leq {1\over \epsilon}\E[|H_n(T)|]\longrightarrow 0.
\end{equation}
Thus, noticing by the definition of $h_n(t)$ together with (22) that
$$
\sup_{t\in\T}\sup_{i\geq 1} \E\left[|\hat{y}_{t}^{n,i}|\right]
\leq \sup_{t\in\T} h_n(t)\leq \Theta^{-1}(\Theta(\sup_{t\in\T}H_n(t))+T)
$$
and
$$
\sup_{t\in\T}\sup_{i\geq 1} |\hat{y}_{t}^{n,i}|=\sup_{t\in\T}\sup_{i\geq 1}\E\left[\left.|\hat{y}_{t}^{n,i}|\right|\F_t\right]
\leq \sup_{t\in\T} h_n(t)\leq \Theta^{-1}(\Theta(\sup_{t\in\T}H_n(t))+T),
\vspace{0.2cm}
$$
from (24) we can deduce that for each $\epsilon>0$,
$$
0\leq \Lim \mathbbm{P}\left(\left\{\sup_{i\geq 1} \sup_{t\in\T}\E\left[|\hat{y}_{t}^{n,i}|\right]\geq \epsilon\right\}\right)\leq
\Lim \mathbbm{P}\left(\left\{ \sup_{t\in\T}|H_n(t)|\geq \Theta^{-1}(\Theta(\epsilon)-T)\right\}\right)=0
$$
and
\begin{equation}
0\leq \Lim \mathbbm{P}\left(\left\{\sup_{i\geq 1}\sup_{t\in\T} |\hat{y}_{t}^{n,i}|\geq \epsilon\right\}\right)\leq
\Lim \mathbbm{P}\left(\left\{ \sup_{t\in\T}|H_n(t)|\geq \Theta^{-1}(\Theta(\epsilon)-T)\right\}\right)=0.
\end{equation}
That is to say,
\begin{equation}
\Lim \sup_{i\geq 1} \sup_{t\in\T} \E\left[|\hat{y}_{t}^{n,i}|\right]=0
\end{equation}
and $\left\{\sup_{i\geq 1}\sup_{t\in\T} |\hat{y}_{t}^{n,i}|\right\}_{n=1}^{+\infty}$
converges to $0$ in probability as $n\To\infty$.\vspace{0.2cm}

Combining (20) and (25) we get that for each $\beta\in (0,1)$,
\begin{equation}
\Lim \sup_{i\geq 1}\E\left[\sup_{t\in\T} |\hat{y}_{t}^{n,i}|^\beta\right]=0.
\end{equation}
It then follows from (19), (26) and (27) that there exists a process $(y_t)_{t\in \T}$ which belongs to the class (D) and the space $\bigcap_{\beta\in (0,1)}\s^\beta$ such that
\begin{equation}
\Lim \sup_{t\in\T}\E\left[|y_t^{n}-y_t|\right]=0
\end{equation}
and for each $\beta\in (0,1)$,
\begin{equation}
\Lim \E\left[\sup_{t\in\T} |y_t^{n}-y_t|^\beta\right]=0.\vspace{0.2cm}
\end{equation}

Furthermore, note that $(\hat y_t^{n,i},\hat z_t^{n,i})_{t\in\T}$ is a solution of the following BSDE:
\begin{equation}
\hat{y}_t^{n,i}=\hat \xi^{n,i}+\int^{T}_{t}\hat g^{n,i}(s,\hat y_s^{n,i})\ \mathrm{d} s-\int^{T}_{t}\hat{z}_s^{n,i}\ \mathrm{d}B_s, \quad t\in\T,
\end{equation}
where for each $n,i\geq 1$ and $y\in\R^k$,
$
\hat g^{n,i}(t,y):=g^{n+i}(t,y+y^n_t)-g^n(t,y^n_t).
$
It follows from assumption (H1) on $g$ together with Remark 1 that $\as$, for each $n,i\geq 1$ and $y\in \R^k$,
\begin{equation}
\begin{array}{lll}
\langle y,\hat g^{n,i}(t,y)\rangle&\leq &\Dis \langle y,g^{n+i}(t,y+y^n_t)-g^{n+i}(t,y^n_t)\rangle
+|y||g^{n+i}(t,y^n_t)-g^n(t,y^n_t)|\vspace{0.1cm}\\
&\leq &\Dis \kappa(|y|^2)+|y||g(t,0)|\mathbbm{1}_{|g(t,0)|>n},
\end{array}
\end{equation}
where the function $\kappa(\cdot)$ is defined in (H1a)$_2$. Thus, by (31), Remark 1 and Lemma 1 we can check that the generator $\hat g^{n,i}$ of BSDE (30) satisfies assumption (A1) with
$
p=\beta,\ \mu=m+2A,\ \nu=0,\ f_t=|g(t,0)|\mathbbm{1}_{|g(t,0)|>n}\ {\rm and}\ \ \varphi_t=\kappa({2A\over m+2A})
$
for each $m\geq 1$ and $\beta\in (0,1)$. It then follows from Proposition 1 with $u=t=0$ and $p=\beta$ that for each $m,n,i\geq 1$, there exist a constant $C_{m,A,\beta,T}>0$ depending only on $m, A,\beta$ and $T$, and a constant $C>0$ such that for each $\beta\in (0,1)$,
\begin{equation}
\begin{array}{lll}
\Dis \E\left[\left(\int_0^T |\hat z_s^{n,i}|^2\ {\rm
d}s\right)^{\beta/2}\right]&\leq &\Dis C_{m,A,\beta,T}\E\left[\sup\limits_{t\in [0,T]}|\hat y_t^{n,i}|^\beta\right]+
C\left(\kappa({2A\over m+2A})\cdot T\right)^{\beta/2}\vspace{0.1cm}\\
&&\Dis +C\E\left[\left(\int^{T}_{0}
|g(t,0)|\mathbbm{1}_{|g(t,0)|>n}\mathrm{d}t\right)^\beta
\right].
\end{array}
\end{equation}
Thus, taking superemum with respect to $i$ and sending first $n\To\infty$ ($m$ being fixed) and then $m\To\infty$ in (32), by virtue of (27), (H2), Lebesgue's dominated convergence theorem and the fact that $\kappa(\cdot)$ is continuous function with $\kappa(0)=0$, we can deduce that for each $\beta\in (0,1)$,
$$
\Lim \sup\limits_{i\geq 1}\E\left[\left(\int_0^T |z^{n+i}_s-z^{n}_s|^2\
{\rm d}s\right)^{\beta/2}\right]=\Lim \sup\limits_{i\geq 1}\E\left[\left(\int_0^T |\hat z^{n,i}_s|^2\ {\rm
d}s\right)^{\beta/2}\right]=0,
$$
which means that there exists a process $(z_t)_{t\in \T}$ which belongs to $\bigcap_{\beta\in (0,1)}\s^\beta$ such that\vspace{0.1cm}
\begin{equation}
\RE\ \beta\in (0,1),\ \ \Lim \E\left[\left(\int_0^T |z^n_s-z_s|^2\
{\rm d}s\right)^{\beta/2}\right]=0.
\end{equation}

Finally, since $\int^{T}_{t}z^n_s\ \mathrm{d} B_s$ converges to $\int^{T}_{t}z_s\ \mathrm{d} B_s$ under the uniform convergence in probability (ucp for short) by (33) and since $y\mapsto g(t,y)$ is continuous and (29) holds true, we can easily check by taking limit in both sides of BSDE (15) under ucp that $(y_t,z_t)_{t\in\T}$ is an $L^1$ solution of BSDE (1). \vspace{0.2cm}\hfill $\Box$

{\bf Remark 3}\ Compared with the existing works, the argument between (18)-(29) involved in Gronwall's inequality, Bihari's ineqality and the relationship between convergence in $L^1$ and in probability seems to be  new. \vspace{0.2cm}

With Proposition 5 in the hand we can prove the main existence result.\vspace{0.2cm}

{\bf Proof of Theorem 2}\ \ Assume that $\xi\in\LT$ and the generator $g$ satisfies assumptions (H1)-(H3). Furthermore, in the case of $\alpha\in [1/2,1)$ in (H3) we also assume that $g$ satisfies (H1) with a function $\rho(\cdot)$ satisfying (14) for some $\bar p>1$. We will use some kind of Picard's iterative procedure. \vspace{0.1cm}

Let us set $(y^0_\cdot,z^0_\cdot):=(0,0)$. Note by (H3) that $g$ has a sublinear growth with respect to $z$. By virtue of H\"{o}lder's inequality and assumptions (H1) and (H2) of $g$, it is not hard to verify that the generator $g(t,y,z_t)$ satisfies (H1) and (H2) for each $(z_t)_{t\in\T}\in\M^1$. Thus, with the help of Proposition 5, we can define the process sequence $\{(y^n_t,z^n_t)_{t\in\T}\}_{n=1}^\infty$ recursively,
\begin{equation}
y^{n}_t=\xi+\int^{T}_{t}g(s,y^{n}_s,z^{n-1}_s)\,\mathrm{d} s-\int^{T}_{t}z^{n}_s\,\mathrm{d} B_s, \quad t\in\T,
\end{equation}
where for each $n\geq 1$, $(y^n_t,z^n_t)_{t\in\T}$ belongs to the space $\s^\beta\times \M^\beta$ for each $\beta\in (0,1)$, and $(y^n_t)_{t\in\T}$ belongs to the class (D).\vspace{0.2cm}

For each $n,i\geq 0$, set $\hat{y}^{n,i}_\cdot:=y^{n+i}_\cdot-y^n_\cdot$ and $\hat z^{n,i}_\cdot:=z^{n+i}_\cdot-z^n_\cdot$. Arguing as in the proof of Theorem 1, we can get, in view of (H1) and (H3), that for each $n,i\geq 1$,
$$|\hat{y}^{n,i}_t|\leq \left(AT+G^{n,i}(t)\right)\cdot e^{AT},\quad t\in\T,$$
where
\begin{equation}
G^{n,i}(t):=2\gamma\E\left[\left.\int^{T}_{0}\left(g_s+
|y^n_s|+|z^{n-1}_s|+|z^{n+i-1}_s|\right)^\alpha\ \mathrm{d}s\right|\F_t\right]\in \s^q
\end{equation}
as soon as $\alpha q<1$ with $q>1$. Hence, for each $n,i\geq 1$, we know that $(\hat y_t^{n,i})_{t\in\T}$ belongs to the space $\s^q$ as soon as $\alpha q<1$ with $q>1$. In the sequel, we will deal with two cases respectively:
$$
(i)\ \ \alpha\in (0,1/2); \quad \quad (ii)\ \ \alpha\in [1/2,1).
$$

{\bf Case $(i)$}: In this case, we can pick $q=2$, then for each $n,i\geq 1$,
\begin{equation}
(\hat y_t^{n,i})_{t\in\T}\in \s^2.
\end{equation}
Note that for each $n,i\geq 1$, $(\hat y_t^{n,i},\hat z_t^{n,i})_{t\in\T}$ is a solution of the following BSDE:
\begin{equation}
\hat{y}_t^{n,i}=\int_t^T \bar{g}^{n,i}(s,\hat{y}_s^{n,i})\ {\rm
d}s-\int_t^T \hat{z}_s^{n,i}{\rm d}B_s,\ \ \ t\in \T,
\end{equation}
where for each $y\in \R^k$, $\bar{g}^{n,i}(s,y):=g(s,y+y_s^n,z^{n+i-1}_s)-
g(s,y_s^n,z^{n-1}_s)$. It follows from (H1), (H3) and Remark 1 that $\as$, for each $y\in\R^k$,
\begin{equation}
\begin{array}{lll}
\langle y,\bar{g}^{n,i}(t,y)\rangle &\leq &\Dis \kappa(|y|^2)+2\gamma|y|\left(g_s+
|y^n_s|+|z^{n-1}_s|+|z^{n+i-1}_s|\right)^\alpha\\
&\leq & \Dis A|y|^2+2\gamma|y|\left(g_s+
|y^n_s|+|z^{n-1}_s|+|z^{n+i-1}_s|\right)^\alpha+A,
\end{array}
\end{equation}
and
\begin{equation}
\langle y,\bar{g}^{n,i}(t,y)\rangle \leq \kappa(|y|^2)+\lam
|y||z^{n+i-1}_t-z^{n-1}_t|=\kappa(|y|^2)+ \lam |y||\hat
z^{n-1,i}_t|,\vspace{0.2cm}
\end{equation}
where $\kappa(\cdot)$ is defined in (H1a)$_2$. Thus, by (35) and (38) we know that the generator $\bar g^{n,i}(t,y)$ of BSDE (37) satisfies assumption (A1) with
$p=2,\ \mu=A,\ \nu=0,\ f_t=2\gamma\left(g_t+
|y^n_t|+|z^{n-1}_t|+|z^{n+i-1}_t|\right)^\alpha\ {\rm and}\ \varphi_t\equiv A$.
Then, in view of (36), it follows from Proposition 1 with $p=2$ that $\hat z^{n,i}_\cdot\in {\rm M}^2$. Consequently, for each $n,i\geq 1$, $(\hat y^{n,i}_t,\hat z^{n,i}_t)_{t\in \T}$ is an $L^2$ solution of BSDE (37).\vspace{0.2cm}

On the other hand, it follows from (39) that for each $n\geq 2$ and $i\geq 1$, the generator $\bar g^{n,i}(t,y)$ of BSDE (37) also satisfies assumption (A2) with
$
p=2,\ \psi(u)=\kappa(u),\ \nu=0,\ f_t=\lam |\hat z^{n-1,i}_t|.
$
Then, by Proposition 3 with $u=0$ and H\"{o}lder's inequality we can deduce that there exists a constant $C>0$ such that for each $t\in [0,T]$,
\begin{equation}
\begin{array}{lll}
&&\Dis \E\left[\sup\limits_{r\in [t,T]}|\hat
y^{n,i}_r|^2\right]+\E\left[\int_t^T |\hat z^{n,i}_s|^2\ {\rm d}s\right]\vspace{0.1cm}\\
&\leq &\Dis e^{C(T-t)}\left\{\int_t^T\kappa\left(\E\left[\sup_{r\in
[s,T]}|\hat y^{n,i}_r|^2 \right]\right){\rm d}s+\lam^2 (T-t)\E\left[\int_t^T |\hat z^{n-1,i}_s|^2\ {\rm d}s\right]\right\}.
\end{array}
\end{equation}
Now, let
$$
\delta T:=\min\left\{{\ln 2\over C}, {1\over 16\lam^2}, {\ln 2\over 2A}\right\}\ \ {\rm and}\ \
T_j:=(T-j\delta T)\vee 0,\ \ \RE\ j=1,2,\cdots\vspace{0.2cm}
$$
Then for each $t\in [T_1,T]$, we have
\begin{equation}
e^{C(T-t)}\leq 2, \ \ \lam^2e^{C(T-t)}(T-t)\leq {1\over 8},\ \ e^{2A(T-t)}\leq 2.
\end{equation}
Combining (40) with (41) yields that for each $n\geq 2$, $i\geq 1$ and $t\in [T_1,T]$,\vspace{0.2cm}
\begin{equation}
\Dis \E\left[\sup\limits_{r\in [t,T]}|\hat
y^{n,i}_r|^2\right]+\E\left[\int_t^T |\hat z^{n,i}_s|^2\ {\rm d}s\right]
\leq \Dis 2\int_t^T\kappa\left(\E\left[\sup_{r\in
[s,T]}|\hat y^{n,i}_r|^2 \right]\right){\rm d}s+{1\over 8}\E\left[\int_t^T |\hat z^{n-1,i}_s|^2\ {\rm d}s\right].\vspace{0.2cm}
\end{equation}

Furthermore, note by Remark 1 that $\kappa(x)\leq A(x+1)$ for each $x\geq 0$. Gronwall's inequality with (42) and (41) yields that for each $n\geq 2$, $i\geq 1$  and $t\in [T_1,T]$,
\begin{equation}
\begin{array}{lll}
\Dis \E\left[\sup\limits_{r\in [t,T]}|
\hat y^{n,i}_r|^2\right]+\E\left[\int_t^T | \hat z^{n,i}_s|^2\ {\rm
d}s\right]&\leq &\Dis \left(2AT+{1\over 8}\E\left[\int_t^T |\hat z^{n-1,i}_s|^2\ {\rm d}s\right]\right)\cdot e^{2A(T-t)}\vspace{0.1cm}\\
&\leq &\Dis 4AT+{1\over 4}\E\left[\int_t^T
|\hat z^{n-1,i}_s|^2\ {\rm d}s\right].
\end{array}
\end{equation}
By picking $n=2$ and $i=m-2$ in (43) we get that for each $t\in [T_1,T]$ and $m\geq 3$,
$$
\begin{array}{lll}
\Dis \E\left[\int_t^T | z^m_s-z^2_s|^2\ {\rm
d}s\right]&\leq &\Dis 4AT+{1\over 4}\E\left[\int_t^T | z^{m-1}_s-z^1_s|^2\ {\rm d}s\right]\vspace{0.1cm}\\
&\leq &\Dis 4AT+{1\over 2}\E\left[\int_0^T |z^2_s- z^1_s|^2\ {\rm d}s\right]+{1\over 2}\E\left[\int_t^T | z^{m-1}_s-z^2_s|^2\ {\rm d}s\right],
\end{array}
$$
from which we can obtain by induction that for each $t\in [T_1,T]$,
\begin{equation}
\sup_{m\geq 3} \E\left[\int_t^T | z^m_s-z^2_s|^2\ {\rm
d}s\right]\leq 8AT+\E\left[\int_0^T |z^2_s-z^1_s|^2\ {\rm d}s\right]<+\infty.
\end{equation}
In addition, for each $n\geq 2$, $i\geq 1$ and $t\in \T$ we have\vspace{0.1cm}
\begin{equation}
\Dis{1\over 4}\E\left[\int_t^T
|\hat z^{n-1,i}_s|^2\ {\rm d}s\right]\leq \Dis {1\over 2}\E\left[\int_t^T
\left(|z^{n-1+i}_s-z^2_s|^2+|z^{n-1}_s-z^2_s|^2\right)\ {\rm d}s\right]\leq \Dis \sup_{m\geq 1} \E\left[\int_t^T | z^m_s-z^2_s|^2\ {\rm d}s\right]\vspace{-0.1cm}.
\end{equation}
Combining (43), (45) and (44) yields that for each $t\in [T_1,T]$,
\begin{equation}
\Dis\sup_{n\geq 2}\sup_{i\geq 1}\left(\E\left[\sup\limits_{r\in [t,T]}|
\hat y^{n,i}_r|^2\right]+\E\left[\int_t^T |\hat z^{n,i}_s|^2\ {\rm d}s\right]\right)\leq \Dis 12AT+\E\left[\int_0^T |z^2_s- z^1_s|^2\ {\rm d}s\right]<+\infty.\vspace{0.1cm}
\end{equation}

Now, in view of (46), by first taking supremum with respect to $i$ and then taking $\limsup$ with respect to $n$ in (42) and finally using Fatou's lemma, the monotonicity and continuity of the function $\kappa(\cdot)$ together with Bihari's inequality, we can deduce the existence of processes $(Y_t,Z_t)_{t\in [T_1,T]}\in {\s}^2(T_1,T;\R^k)\times {\rm M}^2(T_1,T;\R^{k\times d})$ such that
\begin{equation}
\Lim \E\left[\sup_{t\in [T_1,T]}|(y^n_t-y^1_t)-Y_t|^2+\int_{T_1}^T
|(z^n_t-z^1_t)-Z_t|^2\ \mathrm{d}t\right]=0.
\end{equation}
Thus, note that $(y^1_t,z^1_t)_{t\in\T}\in \s^\beta\times \M^\beta$ for each $\beta\in (0,1)$ and $(y^1_t)_{t\in\T}$ belongs to the class (D). By passing to the limit in ucp for BSDE (34), in view of (47), (H2), (H3) and Lebesgue's dominated convergence theorem, we deduce that
$
(y_t,z_t)_{t\in [T_1,T]}:=(Y_t+y^1_t,Z_t+z^1_t)_{t\in [T_1,T]}
$
is an $L^1$ solution to the
BSDE with parameters $(\xi,T,g)$ on $[T_1,T]$.\vspace{0.1cm}

Finally, noticing that the $\delta T>0$ depends only on $\lam$ and $A$, we can find a minimal integer $N\geq 1$ such that $T_N=0$. Thus, we can repeat, in finite steps, the above procedure to obtain an $L^1$ solution to BSDE (1) on $[T_2,T_1]$, $[T_3,T_2]$, $\cdots$, $[0,T_{N-1}]$, and then we find an $L^1$ solution to BSDE (1) on $[0,T]$.\vspace{0.2cm}

{\bf Case $(ii)$}: In this case, we can pick a $q\in (1,\bar p\wedge {1\over \alpha})$, then for each $n,i\geq 1$,
\begin{equation}
(\hat y_t^{n,i})_{t\in\T}\in \s^q.
\end{equation}
Note that $q<\bar p$ and that we also assume that equality (14) holds in this case. It follows from Proposition 1 in \citet{Fan15} that $g$ also satisfies (H1b)$_{\bar p}$ and then (H1b)$_q$. Note further that (37) and (38) is also true. In the same way as that in case $(i)$, it follows from Proposition 1 with $p=q$ that, in view of (48) and (35), $\hat z^{n,i}_t\in {\rm M}^q$. Consequently, for each $n,i\geq 1$, $(\hat y^{n,i}_t,\hat z^{n,i}_t)_{t\in \T}$ is an $L^q$ solution of BSDE (37).\vspace{0.2cm}

Furthermore, it follows from (H1b)$_q$ and (H3) of the generator $g$ together with Remark 1 that $\as$, for each $y\in \R^k$,\vspace{0.2cm}
$$
\left\langle {y\over |y|}\mathbbm{1}_{|y|\neq 0}, \bar{g}^{n,i}(t,y)\right\rangle \leq \varrho^{1\over q}(|y|^q)+\lam |\hat z^{n-1,i}_t|,
$$
where $\varrho(\cdot)$ is defined in (H1b)$_q$. Then, for each $n\geq 2$ and $i\geq 1$, the generator $\bar g^{n,i}(t,y)$ of BSDE (37) satisfies assumption (A3) with
$
p=q,\ \phi(u)=\varrho(u),\ \nu=0,\ f_t=\lam |\hat z^{n-1,i}_t|.
$
Then, by Proposition 4 with $u=0$ and $p=q$ together with H\"{o}lder's inequality we can deduce that there exists a constant $C_q>0$ depending only on $q$ such that for each $n\geq 2$, $i\geq 1$ and $t\in [0,T]$,
\begin{equation}
\begin{array}{lll}
&& \Dis \E\left[\sup\limits_{r\in [t,T]}|\hat
y^{n,i}_r|^q\right]+\E\left[\left(\int_t^T |\hat z^{n,i}_s|^2\ {\rm d}s\right)^{q/2}\right]\\
&\leq &\Dis e^{C_q(T-t)}\left\{\int_t^T\varrho\left(\E\left[\sup_{r\in
[s,T]}|\hat y^{n,i}_r|^q \right]\right){\rm d}s+\lam^q (T-t)^{q/2}\E\left[\left(\int_t^T |\hat z^{n-1,i}_s|^2\ {\rm d}s\right)^{q/2}\right]\right\}.
\end{array}
\end{equation}
Now, let
$$
\overline{\delta T}:=\min\left\{{\ln 2\over C_q}, \left({1\over 16\lam^q}\right)^{2\over q},{\ln 2\over 2A}\right\}\ \ {\rm and}\ \ \bar T_j:=(T-j\overline{\delta T})\vee 0,\ \ \RE\ j=1,2,\cdots\vspace{0.2cm}
$$
Then for each $t\in [\bar T_1,T]$, we have
\begin{equation}
e^{C_q(T-t)}\leq 2, \ \ \lam^q e^{C_q(T-t)}(T-t)^{q/2}\leq {1\over 8},\ \ e^{2A(T-t)}\leq 2.
\end{equation}
Combining (49) with (50) yields that for each $n\geq 2$, $i\geq 1$ and $t\in [\bar T_1,T]$,
\begin{equation}
\begin{array}{ll}
& \Dis \E\left[\sup\limits_{r\in [t,T]}|\hat
y^{n,i}_r|^q\right]+\E\left[\left(\int_t^T |\hat z^{n,i}_s|^2\ {\rm d}s\right)^{q/2}\right]\vspace{0.1cm}\\
\leq &\Dis 2\int_t^T\varrho\left(\E\left[\sup_{r\in
[s,T]}|\hat y^{n,i}_r|^q \right]\right){\rm d}s+{1\over 8}\E\left[\left(\int_t^T |\hat z^{n-1,i}_s|^2\ {\rm d}s\right)^{q/2}\right].
\end{array}
\end{equation}

Furthermore, note by Remark 1 that $\varrho(x)\leq A(x+1)$ for each $x\geq 0$. Gronwall's inequality with (51) and (50) yields that for each $n\geq 2$, $i\geq 1$ and $t\in [\bar T_1,T]$,
\begin{equation}
\begin{array}{lll}
\Dis \E\left[\sup\limits_{r\in [t,T]}|
\hat y^{n,i}_r|^q\right]+\E\left[\left(\int_t^T |\hat z^{n,i}_s|^2\ {\rm d}s\right)^{q/2}\right]
&\leq &\Dis \left(2AT+{1\over 8}\E\left[\left(\int_t^T | \hat z^{n-1,i}_s|^2\ {\rm d}s\right)^{q/2}\right]\right)\cdot e^{2A(T-t)}\\
&\leq &\Dis 4AT+{1\over 4}\E\left[\left(\int_t^T | \hat z^{n-1,i}_s|^2\ {\rm d}s\right)^{q/2}\right],\vspace{-0.2cm}
\end{array}
\end{equation}
from which, in view of $q\in (1,2)$ and the basic inequality
$$\left(\int_t^T(a_s+b_s)^2\ \mathrm{d}s\right)^{q/2}\leq 2\left[\left(\int_t^T a^2_s\ \mathrm{d}s\right)^{q/2}+\left(\int_t^T b^2_s\ \mathrm{d}s\right)^{q/2}\right]$$
for each $a_s,b_s\in L^2([t,T])$, by a similar argument to that in case $(i)$ we can deduce that for each $n\geq 1$ and $t\in [\bar T_1,T]$,
\begin{equation}
\Dis\sup_{n\geq 2}\sup_{i\geq 1}\left(\E\left[\sup\limits_{r\in [t,T]}|
\hat y^{n,i}_r|^q\right]+\E\left[\left(\int_t^T | \hat z^{n,i}_s|^2\ {\rm
d}s\right)^{q/2}\right]\right)\\
\leq \Dis 12AT+\E\left[\left(\int_0^T |z^2_s- z^1_s|^2\ {\rm d}s\right)^{q/2}\right]<+\infty.
\end{equation}

Now, in view of (53), by first taking supremum with respect to $i$ and then taking $\limsup$ with respect to $n$ in (51) and finally using Fatou's lemma, the monotonicity and continuity of the function $\kappa(\cdot)$ together with Bihari's inequality, we can deduce the existence of processes $(Y_t,Z_t)_{t\in [\bar T_1,T]}\in {\s}^q(\bar T_1,T;\R^k)\times {\rm M}^q(\bar T_1,T;\R^{k\times d})$ such that
\begin{equation}
\Lim \E\left[\sup_{t\in [\bar T_1,T]}|(y^n_t-y^1_t)-Y_t|^q+\left(\int_{\bar T_1}^T
|(z^n_t-z^1_t)-Z_t|^2\mathrm{d}t\right)^{q/2}\right]=0.
\end{equation}
Thus, note that $(y^1_t,z^1_t)_{t\in\T}\in \s^\beta\times \M^\beta$ for each $\beta\in (0,1)$ and $(y^1_t)_{t\in\T}$ belongs to the class (D). By passing to the limit in ucp for BSDE (34), in view of (54), (H2), (H3) and Lebesgue's dominated convergence theorem, we deduce that
$
(y_t,z_t)_{t\in [\bar T_1,T]}:=(Y_t+y^1_t,Z_t+z^1_t)_{t\in [\bar T_1,T]}
$
is an $L^1$ solution to the
BSDE with parameters $(\xi,T,g)$ on $[\bar T_1,T]$.\vspace{0.1cm}

Finally, noticing that the positive real number $\overline{\delta T}$ depends only on $q$, $\lam$ and $A$, we can find a minimal integer $\bar N\geq 1$ such that $T_{\bar N}$=0. Thus, we can repeat, in finite steps, the above procedure to obtain an $L^1$ solution to BSDE (1) on $[\bar T_2,\bar T_1]$, $[\bar T_3,\bar T_2]$, $\cdots$, $[0,\bar T_{\bar N-1}]$, and then we find an $L^1$ solution to BSDE (1) on $[0,T]$. The proof of Theorem 2 is finally completed.\vspace{0.2cm} \hfill $\Box$

{\bf Remark 4}\ \ We would like to mention that it is interesting that in the case of $\alpha$ in (H3) values in $[1/2,1)$, the assumption (H1) in Theorem 2 needs to be replaced with the stronger assumption (H1)${}_{\bar p}$. The main reason is to ensure obtaining the key inequality (49). This indicates the difference between Propositions 3 and 4. In addition, we point out that how to divide appropriately the time interval $\T$ is also one of key problems in the proof of Theorem 2.

\section{Existence and Uniqueness of the solution in $\s^1\times{\rm M^1}$ and examples}

In this section, by virtue of Theorems 1 and 2 we will establish an existence and uniqueness result of the solution in the space $\s^1\times {\rm M}^1$ (a new type of $L^1$ solution) for multidimensional BSDEs with generators of one-sided Osgood type. This is the first time to the best of our knowledge. We will also provide two examples in this section to illustrate our theoretical results.\vspace{0.2cm}

{\bf Theorem 3} Assume that the generator $g$ satisfies assumptions (H1)-(H3). In the case when the $\alpha$ defined in (H3) values in $[1/2,1)$, we also assume that there exists a constant $\bar p>1$ such that the function $\rho(\cdot)$ in (H1) satisfies (14). If the following assumption (H4) holds true:\\

{\bf (H4)}\ $\Dis \E\left[\sup\limits_{t\in\T}\left(\E\left[|\xi|
+\left.\int_0^T |g(s,0,0)|\ {\rm d}s\right|\F_t\right]\right)\right]<+\infty.
\vspace{0.4cm}
$

\noindent then BSDE (1) admits a unique solution $(y_t,z_t)_{t\in\T}$ in $\s^1(0,T;\R^k)\times{\rm M}^1(0,T;\R^{k\times d})$.\vspace{0.3cm}

{\bf Proof}\ \ It is clear that the uniqueness follows from Theorem 1 directly. Note that if (H4) holds true, then $\xi\in \LT$. It follows from Theorems 1-2 that BSDE (1) admits a unique $L^1$ solution $(y_t,z_t)_{t\in\T}$, i.e., $(y_t)_{t\in\T}$ belongs to the class (D) and $(y_t,z_t)_{t\in\T}\in \s^\beta\times {\rm M}^\beta$ for each $\beta\in (0,1)$. Hence, in order to complete the proof of Theorem 3, it remains to show, under (H1)-(H4),
$$(y_t,z_t)_{t\in\T}\in \s^1(0,T;\R^k)\times{\rm M}^1(0,T;\R^{k\times d}).$$
In fact, let us fix $n\geq 1$ and denote $\tau_n$ the stopping time
$$
\tau_n:=\inf\left\{t\in\T:\int^t_{0}|z_s|^2\
\mathrm{d}s\geq n\right\}\wedge T.
$$
By Corollary 2.3 in \cite{Bri03} we know that for each $t\in\T$,\vspace{0.1cm}
\begin{equation}
|y_{t\wedge \tau_n}|\leq |y_{\tau_n}|+\int_{t\wedge \tau_n}^{\tau_n}\left\langle {y_s\over |y_s|}\mathbbm{1}_{|y_s|\neq 0}, g(s,y_s,z_s)\right\rangle{\rm d}s-\int_{t\wedge \tau_n}^{\tau_n}\left\langle {y_s\over |y_s|}\mathbbm{1}_{|y_s|\neq 0},z_s {\rm d}B_s\right\rangle.\vspace{0.1cm}
\end{equation}
And, it follows from (H1) and (H3) that $\ass$,
\begin{equation}
\left\langle {y_s\over |y_s|}\mathbbm{1}_{|y_s|\neq 0}, g(s,y_s,z_s)\right\rangle\leq \rho(|y_s|)+|g(s,0,0)|+\gamma(g_s+|z_s|)^{\alpha}.
\end{equation}
Thus, in view of (55) and (56), using a similar argument to that in the proof of Theorem 1 we can get that for each $t\in \T$,
\begin{equation}
|y_t|\leq (AT+\bar G(t))\cdot e^{AT},
\end{equation}
where
$$\bar G(t):=\E\left[|\xi|+\left.\int_0^T|g(s,0,0)|\ \mathrm{d}s+\gamma\int_0^T(g_s+|z_s|)^{\alpha}\
\mathrm{d}s\right|\F_t\right],\ \ t\in \T.\vspace{0.2cm}
$$
Furthermore, it follows from (H4) and a similar argument to obtain (5) that $\bar G(\cdot)\in \s^1(0,T;\R^k)$ and then, in view of (57), $(y_t)_{t\in\T}\in \s^1(0,T;\R^k)$. Finally, note by (H1) and (H3) together with Remark 1 that $\as$, for each $(y,z)\in \R^k\times\R^{k\times d}$,
$$
\langle y, g(t,y,z) \rangle\leq \bar\kappa(|y|^2)+\lambda |y||z|+|y||g(t,0,0)|\leq A|y|^2+\lambda |y||z|+|y||g(t,0,0)|+A,
$$
where the function $\bar\kappa(\cdot)$ is defined in (H1a)$_2$. It follows from Proposition 1 with $p=1$ and $u=t=0$ that $(z_t)_{t\in\T}\in {\rm M}^1(0,T;\R^{k\times d})$. Theorem 3 is then proved.\vspace{0.2cm} \hfill $\Box$

By Theorems 1-3 and Remarks 1 and 2 the following corollary is immediate.\vspace{0.1cm}

{\bf Corollary 1} Assume that the generator $g$ satisfies assumption (H1b)$_p$ for some $p>1$, (H2) and (H3). Then, for each $\xi\in\LT$, BSDE (1) admits a unique $L^1$ solution $(y_t,z_t)_{t\in\T}$. Furthermore, if (H4) holds true, then
$(y_t,z_t)_{t\in\T}\in \s^1(0,T;\R^k)\times {\rm M}^1(0,T;\R^{k\times d})$.\vspace{0.2cm}

{\bf Remark 5} Note that if the generator $g$ satisfies the monotonicity condition used in \citet{Bri03}, then it must satisfy (H1b)$_p$ for all $p>1$. Theorems 1-3 of this paper generalize Theorems 6.2 and 6.3 in \cite{Bri03}.\vspace{0.2cm}

{\bf Example 1} Let $k=1$ and
$$
g(\omega,t,y,z)=h(|y|)-e^{|B_t(\omega)|\cdot y}+(e^{-y}\wedge 1)\cdot \sin |z|+{1\over \sqrt{t}}\mathbbm{1}_{t>0},
$$
where
$$
h(x)=\left\{
\begin{array}{lll}
-x|\ln x|& ,&0<x\leq \delta;\\
h'(\delta-0)(x-\delta)+h(\delta)& ,&x> \delta;\\
0& ,&{\rm other\ cases}
\end{array}\right.\vspace{0.1cm}
$$
with $\delta>0$ small enough.

It is not hard to check that $g$ satisfies assumptions (H2) and (H3) with $\lam =1$ and any $\alpha\in (0, 1/2)$. Furthermore, note that $e^{-\beta y}$ is decreasing in $y$ for each $\beta\geq 0$, $h(\cdot)$ is concave and sub-additive and then the following inequality holds: $\as$,
$$
\RE\ y_1,y_2,z,\ \
\left\langle {y_1-y_2\over |y_1-y_2|}\mathbbm{1}_{|y_1-y_2|\neq 0},
g(\omega,t,y_1,z)-g(\omega,t,y_2,z)\right\rangle\leq
h(|y_1-y_2|)
$$
with $\int_{0^+} {{\rm d}u\over h(u)}=+\infty$.
It follows that $g$ also satisfies assumption (H1). Then, by Theorems 1-2 we know that for each $\xi\in\LT$, the BSDE with the parameters $(\xi,T,g)$ admits a unique $L^1$ solution $(y_t,z_t)_{t\in\T}$. Moreover, by Theorem 3 we also know that if (H4) holds true for $\xi$ and $g(t,0,0)$, then $$(y_t,z_t)_{t\in\T}\in \s^1(0,T;\R^k)\times {\rm M}^1(0,T;\R^{k\times d}).$$

{\bf Example 2}\ Let $y=(y_1,\cdots,y_k)$ and
$g(\omega,t,y,z)=(g_1(\omega,t,y,z),\cdots,g_k(\omega,t,y,z))$, where for each $i=1,\cdots,k$,
$$
g_i(\omega,t,y,z):=e^{-y_i}+\bar h(|y|)+\left(|z|^2\wedge |z|^{2/3}\right)+|B_t(\omega)|,
$$
with
$$
\bar h(x)=\left\{
\begin{array}{lll}
-x|\ln x|^{1/p}& ,&0<x\leq \delta;\\
\bar h'(\delta-0)(x-\delta)+\bar h(\delta)& ,&x> \delta;\\
0& ,&{\rm other\ cases}
\end{array}\right.\vspace{0.2cm}
$$
with $\delta>0$ small enough and $p>1$.

In the same way as in Example 1, we can check that this generator $g$ satisfies assumptions (H1b)$_p$ with function $\bar h(\cdot)$, (H2) and (H3) with $\lam =1$ and $\alpha=2/3$. It then follows from Corollary 1 that for each $\xi\in\LT$, the BSDE with the parameters $(\xi,T,g)$ admits a unique $L^1$ solution $(y_t,z_t)_{t\in\T}$. And, if (H4) holds true for $\xi$ and $g(t,0,0)$, then $(y_t,z_t)_{t\in\T}\in \s^1(0,T;\R^k)\times {\rm M}^1(0,T;\R^{k\times d}).$

\section{Stability of the $L^1$ solutions and the solutions in $\s^1\times {\rm M}^1$}

In this section, enlightened by the proof of Proposition 5 and Theorems 2-3, we shall put forward and prove the stability theorems of the $L^1$ solutions and the solutions in the space $\s^1\times {\rm M}^1$ for multidimensional BSDEs with generators of one-sided Osgood type. To the best of our knowledge, this is the first time for the $L^1$ solution of multidimensional BSDEs.\vspace{0.2cm}

In the sequel, for each $m\in \N$, let $\xi^m\in \LT$ and let  $(y_t^m, z_t^m)_{t\in [0,T]}$ be
an $L^1$ solution of
the following BSDEs depending on parameter $m$:
\begin{equation}
y_t^m=\xi^m+\int_t^T g^m(s,y_s^m,z_s^m)\ {\rm d}s -\int_t^T z_s^m{\rm d}B_s,\quad t\in [0,T].
\end{equation}
Furthermore,  we introduce  the following assumptions:\vspace{0.2cm}

{\bf (B1)}\ All $g^m$ satisfy assumptions (H1)-(H3)
 with the same parameters $\rho(\cdot)$, $\lambda$, $\gamma$, $g_t$ and $\alpha$. Furthermore, in the case of $\alpha\in [1/2,1)$ we assume that $g$ satisfies (H1) with a function $\rho(\cdot)$ satisfying (14) for some $\bar p>1$. \vspace{0.2cm}

{\bf (B2)}\ There exists a nonnegative real number sequence
$\{a_m\}_{m=1}^{+\infty}$ satisfying $\lim\limits_{m\To \infty}a_m=0$ such that
$\as$,  for each $m\geq 1$,
\begin{equation}
\RE\ (y,z)\in \R^k\times \R^{k\times d},\ \ |g^m(\omega,t,y,z)-g^0(\omega,t,y,z)|\leq a_m.
\end{equation}
And,
\begin{equation}
\lim\limits_{m\To \infty}\E\left[|\xi^{m}-\xi^{0}|\right]=0.
\vspace{0.3cm}
\end{equation}

The following Theorem 4 is the stability theorem of $L^1$ solutions.\vspace{0.1cm}

{\bf Theorem 4}\ Under assumptions (B1) and (B2), we have
\begin{equation}
\lim\limits_{m\To \infty}\sup\limits_{t\in [0,T]}\E\left[|y_t^m-y_t^0|\right]=0,
\end{equation}
and for each $\beta\in (0,1)$,
\begin{equation}
\lim\limits_{m\To \infty}\E\left[\sup\limits_{t\in [0,T]} |y_t^m-y_t^0|^\beta+\left(\int_0^T |z_s^m -z_s^0|^2 {\rm d}s\right)^{\beta/2} \right]=0.\vspace{0.1cm}
\end{equation}

{\bf Proof}\ \ For each $ m\in \N$, set $(y^{m,0}_\cdot,z^{m,0}_\cdot):=(0,0)$ and, similar to the beginning part of the proof of Theorem 2, define recursively the process sequence $\{(y^{m,n}_\cdot,z^{m,n}_\cdot)\}_{n=1}^{+\infty}$ by the $L^1$ solutions of the following BSDEs
\begin{equation}
y_t^{m,n}=\xi^m+\int_t^T g^{m,n}(s,y_s^{m,n},z_s^{m,n-1})\ {\rm d}s -\int_t^T z_s^{m,n}{\rm d}B_s,\ \  t\in [0,T],\vspace{0.2cm}
\end{equation}
where $(y^{m,n}_t,z^{m,n}_t)_{t\in [0,T]}\in \s^\beta\times \M^\beta$ for each $\beta\in (0,1)$ and $(y^{m,n}_t)_{t\in\T}$ belongs to the class (D) for each $m,n\in \N$.\vspace{0.2cm}

In the sequel, note that $\as$,
\begin{equation}
|y_t^m -y_t^0| \leq |y_t^m -y_t^{m,n}| + |y_t^{m,n} -y_t^{0,n}|+|y_t^{0,n} -y_t^0|\leq \Dis 2\sup_{m\geq 0}|y_t^{m,n}-y_t^m|+|y_t^{m,n} -y_t^{0,n}|
\end{equation}
and
\begin{equation}
|z_t^m -z_t^0|\leq |z_t^m -z_t^{m,n}| + |z_t^{m,n} -z_t^{0,n}|+|z_t^{0,n} -z_t^0|\leq \Dis 2\sup_{m\geq 0}|z_t^{m,n}-z_t^m|+|z_t^{m,n} -z_t^{0,n}|.
\end{equation}
We will estimate, respectively, every term of the right hand in (64) and (65).\vspace{0.2cm}

Firstly, the following Proposition 6 gives the  estimates with respect to the second term of the right hand in (64) and (65).\vspace{0.1cm}

{\bf Proposition 6}\ \ For each $n\geq 1$, we have
\begin{equation}
\lim\limits_{m\To \infty}\sup\limits_{t\in [0,T]}\E\left[|y_t^{m,n}-y_t^{0,n}|\right]=0,
\end{equation}
and for each $\beta\in (0,1)$,
\begin{equation}
\lim\limits_{m\To \infty}\E\left[\sup\limits_{t\in [0,T]} |y_t^{m,n}-y_t^{0,n}|^\beta+\left(\int_0^T |z_s^{m,n}-z_s^{0,n}|^2\  {\rm d}s\right)^{\beta/2} \right]=0.
\end{equation}

{\bf Proof}\ \ We first consider the case of $n=1$. Let us fix $k,m\geq 1$ and denote $\tau^m_k$ the stopping time
$$
\tau^m_k:=\inf\left\{t\in\T:\int^t_{0}(|z^{m,1}_s|^2
+|z^{0,1}_s|^2)\ \mathrm{d}s\geq k\right\}\wedge T.
$$
Corollary 2.3 in \cite{Bri03} leads to the following inequality
\begin{equation}
\begin{array}{lll}
\Dis|\hat{y}^{m,1}_{t\wedge\tau^m_k}|&\leq &\Dis |\hat{y}^{m,1}_{\tau^m_k}|+\int^{\tau^m_k}_{t\wedge\tau^m_k}
\left\langle {\hat{y}^{m,1}_s\over |\hat{y}^{m,1}_s| }\mathbbm{1}_{|\hat{y}^{m,1}_s|\neq 0},g^m(s,y^{m,1}_s,0)-g^0(s,y^{0,1}_s,0)\right\rangle\
\mathrm{d}s\vspace{0.1cm}\\
&&\Dis -\int^{\tau^m_k}_{t\wedge\tau^m_k}\left\langle {\hat{y}^{m,1}_s\over |\hat{y}^{m,1}_s|}\mathbbm{1}_{|\hat{y}^{m,1}_s|\neq 0},\hat{z}^{m,1}_s\ \mathrm{d} B_s\right\rangle,\ \ t\in\T,
\end{array}
\end{equation}
where and hereafter
$$\hat{y}^{m,1}_\cdot:=y^{m,1}_\cdot-y^{0,1}_\cdot\ \ {\rm and}\ \ \hat{z}^{m,1}_\cdot:=z^{m,1}_\cdot-z^{0,1}_\cdot.$$
And, it follows from assumption (H1) of $g^m$  and (59) that $\ass$,
\begin{equation}
\begin{array}{ll}
&\Dis \left\langle{\hat{y}^{m,1}_s\over |\hat{y}^{m,1}_s|} \mathbbm{1}_{|\hat{y}^{m,1}_s|\neq 0},g^m(s,y^{m,1}_s,0)-g^0(s,y^{0,1}_s,0)
\right\rangle\vspace{0.1cm}\\
\leq &\Dis \left\langle{\hat{y}^{m,1}_s\over |\hat{y}^{m,1}_s|} \mathbbm{1}_{|\hat{y}^{m,1}_s|\neq 0},g^m(s,y^{m,1}_s,0)-g^m(s,y^{0,1}_s,0)
\right\rangle+|g^m(s,y^{0,1}_s,0)-g^0(s,y^{0,1}_s,0)|
\vspace{0.1cm}\\
\leq &\Dis \rho(|\hat{y}^{m,1}_s|)+a_m.
\end{array}
\end{equation}
Then, combining (68) with (69) we can deduce that for each $k,m\geq 1$,
\begin{equation}
|\hat{y}^{m,1}_{t\wedge\tau^m_k}|\leq a_mT+\E\left[
|\hat{y}^{m,1}_{\tau^m_k}|
+\left.\int^{\tau^m_k}_{t\wedge\tau^m_k}
\rho(|\hat{y}^{m,1}_s|)\ \mathrm{d}s\right|\F_t\right],\ \ t\in\T.
\end{equation}
Since $\hat y^{m,1}_\cdot$ belongs to the class (D), and $\rho(\cdot)$ increases at most linearly, we can send $k$ to $+\infty$ in (70) and use Lebesgue's dominated convergence theorem, in view of $\tau^m_k\To T$ as $k\To \infty$ and $\hat y^{m,1}_T=\xi^m-\xi^0$, to get that for each $m\geq 1$,
\begin{equation}
|\hat{y}^{m,1}_t| \leq H_m(t) +\E\left[\left.\int_t^T\rho(|\hat{y}^{m,1}_s|)\ \mathrm{d}s\right|\F_t\right],\ \ t\in \T,
\end{equation}
where
$$H_m(t):=a_mT+\E\left[\left.|\xi^m-\xi^0|
\right|\F_t\right].\vspace{0.2cm}$$
In the sequel, note by Lemma 6.1 in \cite{Bri03} and assumption (B2) that
\begin{equation}
\RE\ \beta\in (0,1),\ \ \sup_{m\geq 1}\E\left[\sup_{t\in\T}|H_m(t)|^\beta
\right]\leq {1\over 1-\beta}\sup_{m\geq 1} \left(\E\left[H_m(T)\right]\right)^\beta<+\infty
\end{equation}
and
\begin{equation}
\lim\limits_{m\To \infty} \E\left[|H_m(T)|\right]=0.\vspace{0.2cm}
\end{equation}
Arguing as that from (18) to (29) we can deduce that for each $\beta\in (0,1)$,
\begin{equation}
\lim\limits_{m\To \infty}\left(\sup\limits_{t\in [0,T]}\E\left[|\hat y_t^{m,1}|\right]+
\E\left[\sup\limits_{t\in [0,T]} |\hat y_t^{m,1}|^\beta\right]\right)=0.
\end{equation}
Furthermore, note that for each $m\geq 1$, $(\hat y_t^{m,1},\hat z_t^{m,1})_{t\in\T}$ is an $L^1$ solution of the following BSDE:
\begin{equation}
\hat{y}_t^{m,1}=\xi^{m}-\xi^0+\int^{T}_{t}\hat g^{m,1}(s,\hat y_s^{m,1})\ \mathrm{d} s-\int^{T}_{t}\hat{z}_s^{m,1}\mathrm{d}B_s, \quad t\in\T,
\end{equation}
where for each $y\in\R^k$,
$
\hat g^{m,1}(t,y):=g^{m}(t,y+y^{0,1}_t,0)-g^0(t,y^{0,1}_t,0).
$
It follows from assumption (H1) on $g^m$ together with (59) that $\as$, for each $m\geq 1$ and $y\in \R^k$,
\begin{equation}
\begin{array}{lll}
\langle y,\hat g^{m,1}(t,y)\rangle&\leq &\Dis \langle y,g^{m}(t,y+y^{0,1}_t,0)-g^{m}(t,y^{0,1}_t,0)\rangle +|y||g^{m}(t,y^{0,1}_t,0)-g^0(t,y^{0,1}_t,0)|\\
&\leq &\Dis \kappa(|y|^2)+|y|a_m,
\end{array}
\end{equation}
where the function $\kappa(\cdot)$ is defined in (H1a)$_2$. Thus, in view of (74)-(76), Proposition 1, (B2) and the assumption of $\kappa(\cdot)$, a similar argument to that from (30) to (33) yields that for each $\beta\in (0,1)$,
\begin{equation}
\lim\limits_{m\To \infty} \E\left[\left(\int_0^T |\hat z^{m,1}_s|^2\
{\rm d}s\right)^{\beta/2}\right]=0.
\end{equation}
From (74) and (77), we know that (66) and (67) hold true for $n=1$.\vspace{0.2cm}

Now, let us fix arbitrarily a $n\geq 2$ and assume that (66) and (67) hold true for $n-1$. In the sequel, we will prove that they also hold for $n$. For each $m\geq 1$, define
$$\hat{y}^{m,n}_\cdot:=y^{m,n}_\cdot-y^{0,n}_\cdot\ \ {\rm and}\ \ \hat{z}^{m,n}_\cdot:=z^{m,n}_\cdot-z^{0,n}_\cdot.$$
Then, $(\hat y_t^{m,n},\hat z_t^{m,n})_{t\in\T}$ is an $L^1$ solution of the following BSDE:
\begin{equation}
\hat{y}_t^{m,n}=\xi^{m}-\xi^0+\int^{T}_{t}\hat g^{m,n}(s,\hat y_s^{m,n})\ \mathrm{d} s-\int^{T}_{t}\hat{z}_s^{m,n}\ \mathrm{d}B_s, \quad t\in\T,
\end{equation}
where for each $y\in\R^k$,
$
\hat g^{m,n}(t,y):=g^{m}(t,y+y^{0,n}_t,z^{m,n-1}_t)
-g^0(t,y^{0,n}_t,z^{0,n-1}_t).
$
It follows from assumption (H1) of $g^m$ together with Remark 1 and (59) that $\as$, for each $m\geq 1$ and $y\in \R^k$,
\begin{equation}
\begin{array}{lll}
\langle y,\hat g^{m,n}(t,y)\rangle&\leq &\Dis \langle y,g^{m}(t,y+y^{0,n}_t,z^{m,n-1}_t)-
g^{m}(t,y^{0,n}_t,z^{m,n-1}_t)\rangle\vspace{0.1cm}\\
&&\Dis +|y||g^{m}(t,y^{0,n}_t,z^{m,n-1}_t)
-g^m(t,y^{0,n}_t,z^{0,n-1}_t)|\vspace{0.1cm}\\
&&\Dis +|y||g^{m}(t,y^{0,n}_t,z^{0,n-1}_t)
-g^0(t,y^{0,n}_t,z^{0,n-1}_t)|\vspace{0.1cm}\\
&\leq &\Dis \kappa(|y|^2)+|y|\cdot(\triangle_t^{m,n}+a_m),
\end{array}
\end{equation}
where the function $\kappa(\cdot)$ is defined in (H1a)$_2$, and, in view of (H3),
\begin{equation}
\triangle_t^{m,n}:=|g^{m}(t,y^{0,n}_t,z^{m,n-1}_t)
-g^m(t,y^{0,n}_t,z^{0,n-1}_t)|\leq \lambda |z^{m,n-1}_t-z^{0,n-1}_t|
\end{equation}
as well as
\begin{equation}
\triangle_t^{m,n}\leq 2\gamma (g_t+|y^{0,n}_t|+|z^{m,n-1}_t|+|z^{0,n-1}_t|)^\alpha.
\end{equation}
Since (67) holds true for $n-1$, from (80) we know that the sequence of random variables $\{\int_0^T \triangle_t^{m,n}\ {\rm d}t\}_{m=1}^{+\infty}$ converges in probability to $0$ as $m\To \infty$, and from (81) and H\"{o}lder's inequality that for each $q>1$ satisfying $\alpha q<1$,
\begin{equation}
\sup_{m\geq 1}\E\left[\left(\int_0^T \triangle_t^{m,n}\ {\rm d}t\right)^q\right]\leq K_n\left(1+\sup_{m\geq 1}\E\left[\left(\int_0^T|z^{m,n-1}_t|^2 \ {\rm d}t\right)^{\alpha q\over 2}\right] \right)<+\infty,
\end{equation}
where $K_n>0$ is a constant independent of $m$. Hence, for each $q'>1$ satisfying $\alpha q'<1$, we have
\begin{equation}
\lim\limits_{m\To \infty} \E\left[\left(\int_0^T\triangle_t^{m,n} \
{\rm d}t\right)^{q'}\right]=0.\vspace{0.1cm}
\end{equation}
On the other hand, note from assumption (H1) of $g^m$  and (79) that $\as$,
\begin{equation}
\left\langle{\hat{y}^{m,n}_t\over |\hat{y}^{m,n}_t|} \mathbbm{1}_{|\hat{y}^{m,n}_t|\neq 0},g^m(t,y^{m,n}_t,z^{m,n-1}_t)
-g^0(t,y^{0,n}_t,z^{0,n-1}_t)
\right\rangle\leq \rho(|\hat{y}^{m,n}_t|)+\triangle_t^{m,n}+a_m.
\end{equation}
Arguing as that from (68) to (71), in view of (82), (83) and (84), we can obtain that for each $m\geq 1$,
\begin{equation}
|\hat{y}^{m,n}_t| \leq H^n_m(t) +\E\left[\left.\int_t^T\rho(|\hat{y}^{m,n}_s|)\ \mathrm{d}s\right|\F_t\right],\ \ t\in \T,
\end{equation}
where
$$H^n_m(t):=a_mT+\E\left[\left.|\xi^m-\xi^0|
+\int_0^T\triangle_s^{m,n} \
{\rm d}s\right|\F_t\right].\vspace{0.1cm}$$
Furthermore, in view of (72), (73), (82), (83) and (85), a similar argument to that from (18) to (29) yields that for each $\beta\in (0,1)$,
\begin{equation}
\lim\limits_{m\To \infty}\left(\sup\limits_{t\in [0,T]}\E\left[|\hat y_t^{m,n}|\right]+
\E\left[\sup\limits_{t\in [0,T]} |\hat y_t^{m,n}|^\beta\right]\right)=0.
\end{equation}
Finally, in view of (78), (79), (83), Proposition 1, (B2) and the assumption of $\kappa(\cdot)$, using a similar argument to that from (30) to (33) yields that for each $\beta\in (0,1)$,
\begin{equation}
\lim\limits_{m\To \infty} \E\left[\left(\int_0^T |\hat z^{m,n}_s|^2\
{\rm d}s\right)^{\beta/2}\right]=0.
\end{equation}
In view of (86) and (87), we have proved that (66) and (67) hold also true for $n$. Thus, by induction the proof of Proposition 6 is completed.\hfill $\Box$\vspace{0.2cm}

Next, we turn to the estimates with respect to the first term of the right hand in (64) and (65).\vspace{0.2cm}

{\bf Proposition 7}\ \ In the case when the $\alpha$ defined in (H3) values in $(0,1/2)$, there exists a positive real number $\delta T>0$ depending only on $\lam$ and $A$ such that for each $t\in [T_1, T]$ with $T_1:=(T-\delta T)\vee 0$,
\begin{equation}
\Lim \sup_{m\geq 0}\E\left[\sup_{r\in [t,T]}|y_r^{m,n}-y_r^m|^2+\int_t^T |z_s^{m,n}-z_s^m|^2\ {\rm d}s\right]=0.
\end{equation}
In the case when the $\alpha$ defined in (H3) values in $[1/2,1)$, for each $q\in (1,\bar p\wedge {1\over \alpha})$, there exists a positive real number $\overline{\delta T}>0$ depending only on $q$, $\lam$ and $A$ such that for each $t\in [\bar T_1, T]$ with $\bar T_1:=(T-\overline{\delta T})\vee 0$,
\begin{equation}
\Lim \sup_{m\geq 0}\E\left[\sup_{r\in [t,T]}|y_r^{m,n}-y_r^m|^q+\left(\int_t^T |z_s^{m,n}-z_s^m|^2\ {\rm d}s\right)^{q/2}\right]=0.\vspace{0.2cm}
\end{equation}

{\bf Proof}\ \ We only prove the case when the $\alpha$ in (H3) values in $(0,1/2)$. In view of the proof of case $(ii)$ in Theorem 2, another case can be proved in the same way.\vspace{0.2cm}

Now, we assume that $\alpha\in (0,1/2)$ and set $\bar{y}^{m,n}_\cdot:=y^{m,n}_\cdot-y^m_\cdot$ and $\bar z^{m,n}_\cdot:=z^{m,n}_\cdot-z^m_\cdot$ for each $m,n\geq 0$. Note that for each $m\geq 0$ and $n\geq 1$, $(\bar y_t^{m,n},\bar z_t^{m,n})_{t\in\T}$ is a solution of the following BSDE:
\begin{equation}
\bar{y}_t^{m,n}=\int_t^T \bar{g}^{m,n}(s,\bar{y}_s^{m,n})\ {\rm
d}s-\int_t^T \bar{z}_s^{m,n}{\rm d}B_s,\ \ \ t\in \T,
\end{equation}
where for each $y\in \R^k$,
$\bar{g}^{m,n}(s,y):=g^m(s,y+y_s^m,z^{m,n-1}_s)-
g^m(s,y_s^m,z^m_s).$
It follows from (H1) of $g^m$ and Remark 1 that $\as$, for each $y\in\R^k$,\vspace{-0.2cm}
\begin{equation}
\begin{array}{lll}
\langle y,\bar{g}^{m,n}(t,y)\rangle &\leq &\langle y,g^m(t,y+y_t^m,z^{m,n-1}_t)-g^m(t,y_t^m,z^{m,n-1}_t)
\rangle\\
&& \Dis +|y||g^m(t,y_t^m,z^{m,n-1}_t)-g^m(t,y_t^m,z^m_t)|\\
&\leq & \Dis \kappa(|y|^2)+|y|\bar\triangle^{m,n}_t
\end{array}
\end{equation}
where $\kappa(\cdot)$ is defined in (H1a)$_2$, and in view of (H3) of $g^m$,
\begin{equation}
\bar\triangle^{m,n}_t:=g^m(t,y_t^m,z^{m,n-1}_t)-g^m(t,y_t^m,z^m_t)\leq \lam |z^{m,n-1}_t-z^m_t|=\lam |\bar z^{m,n-1}_t|
\end{equation}
as well as
\begin{equation}
\bar\triangle^{m,n}_t\leq 2\gamma \left(g_t+
|y^m_t|+|z^m_t|+|z^{m,n-1}_t|\right)^\alpha.\vspace{0.2cm}
\end{equation}
In view of $\alpha\in (0,1/2)$ and (93), using H\"{o}lder's inequality, Jensen's inequality and Doob's inequality yields that for each $m\geq 0$ and $n\geq 1$,
\begin{equation}
\E\left[\left(\int_0^T \bar\triangle^{m,n}_t\ {\rm d}t\right)^2\right]<+\infty\ \ {\rm and\ then}\ \
\E\left[\left.\int_0^T \bar\triangle^{m,n}_s\ {\rm d}s\right|\F_t\right]\in \s^2.
\end{equation}
On the other hand, note from assumption (H1) of $g^m$  and (91) that $\as$,
\begin{equation}
\left\langle{\bar{y}^{m,n}_t\over |\bar{y}^{m,n}_t|} \mathbbm{1}_{|\bar{y}^{m,n}_t|\neq 0},g^m(t,y^{m,n}_t,z^{m,n-1}_t)
-g^m(t,y^m_t,z^m_t)
\right\rangle\leq \rho(|\bar{y}^{m,n}_t|)+\bar\triangle_t^{m,n}.
\end{equation}
Thus, in view of (95) and (94), arguing as that from (2) to (7), by virtue of Gronwall's inequality we can obtain that for each $m\geq 0$ and $n\geq 1$,
\begin{equation}
(\bar y_t^{m,n})_{t\in\T}\in \s^2.
\end{equation}
Furthermore, in view of (96), (91), (94) and Remark 1, it follows from Proposition 1 that $(\bar z_t^{m,n})_{t\in\T}\in {\rm M}^2$. Thus, $(\bar y_t^{m,n},\bar z_t^{m,n})_{t\in\T}$ is an $L^2$ solution of BSDE (90) for each $m\geq 0$ and $n\geq 1$.\vspace{0.1cm}

In the sequel, in view of (91) and (92), by Proposition 3 and H\"{o}lder's inequality we can deduce the existence of a constant $C>0$ such that for each $m\geq 0$, $n\geq 2$ and $t\in [0,T]$,
\begin{equation}
\hspace*{-0.8cm}\begin{array}{lll}
&&\Dis \E\left[\sup\limits_{r\in [t,T]}|\bar
y^{m,n}_r|^2\right]+\E\left[\int_t^T |\bar z^{m,n}_s|^2\ {\rm d}s\right]\vspace{0.1cm}\\
&\leq &\Dis e^{C(T-t)}\left\{\int_t^T\kappa\left(\E\left[\sup_{r\in
[s,T]}|\bar y^{m,n}_r|^2 \right]\right){\rm d}s+\lam^2 (T-t)\E\left[\int_t^T |\bar z^{m,n-1}_s|^2\ {\rm d}s\right]\right\}.
\end{array}
\end{equation}
Now, let
$$
\delta T:=\min\left\{{\ln 2\over C}, {1\over 16\lam^2}, {\ln 2\over 2A}\right\}\ \ {\rm and}\ \
T_1:=(T-\delta T)\vee 0.\vspace{0.2cm}
$$
Then for each $t\in [T_1,T]$, we have
\begin{equation}
e^{C(T-t)}\leq 2, \ \ \lam^2e^{C(T-t)}(T-t)\leq {1\over 8},\ \ e^{2A(T-t)}\leq 2.
\end{equation}
Combining (97) with (98) yields that for each $m\geq 0$, $n\geq 2$ and $t\in [T_1,T]$,\vspace{0.2cm}
\begin{equation}
\Dis \E\left[\sup\limits_{r\in [t,T]}|\bar
y^{m,n}_r|^2\right]+\E\left[\int_t^T |\bar z^{m,n}_s|^2\ {\rm d}s\right]
\leq \Dis 2\int_t^T\kappa\left(\E\left[\sup_{r\in
[s,T]}|\bar y^{m,n}_r|^2 \right]\right){\rm d}s+{1\over 8}\E\left[\int_t^T |\bar z^{m,n-1}_s|^2\ {\rm d}s\right].
\end{equation}
Furthermore, note by Remark 1 that $\kappa(x)\leq A(x+1)$ for each $x\geq 0$. Gronwall's inequality with (99) and (98) yields that for each $m\geq 0$, $n\geq 2$ and $t\in [T_1,T]$,
\begin{equation}
\begin{array}{lll}
\Dis \E\left[\sup\limits_{r\in [t,T]}|
\bar y^{m,n}_r|^2\right]+\E\left[\int_t^T | \bar z^{m,n}_s|^2\ {\rm
d}s\right]&\leq &\Dis \left(2AT+{1\over 8}\E\left[\int_t^T |\bar z^{m,n-1}_s|^2\ {\rm d}s\right]\right)\cdot e^{2A(T-t)}\vspace{0.1cm}\\
&\leq &\Dis 4AT+{1\over 4}\E\left[\int_t^T
|\bar z^{m,n-1}_s|^2\ {\rm d}s\right].
\end{array}
\end{equation}
By the above inequality (100) together with the inequality
$$|\bar z^{m,n-1}_s|^2\leq 2(|z^{m,n}_s-z_s^m|^2+|z^{m,n}_s-z_s^{m,n-1}|^2)
=2(|\bar z^{m,n}_s|^2+|\tilde{z}^{m,n}_s|^2),$$
we can get that for each $m\geq 0$, $n\geq 2$ and $t\in [T_1,T]$,
\begin{equation}
\E\left[\sup\limits_{r\in [t,T]}|
\bar y^{m,n}_r|^2\right]+{1\over 2}\E\left[\int_t^T | \bar z^{m,n}_s|^2\ {\rm d}s\right]\leq 4AT+{1\over 2}\E\left[\int_t^T
|\tilde{z}^{m,n}_s|^2\ {\rm d}s\right],
\end{equation}
where and hereafter, for each $m\geq 0$ and $n\geq 1$ we define
$$\tilde{y}^{m,n}_\cdot:=y^{m,n}_\cdot-y^{m,n-1}_\cdot\ \ {\rm and}\ \ \tilde{z}^{m,n}_\cdot:=z^{m,n}_\cdot-z^{m,n-1}_\cdot.
\vspace{-0.1cm}$$

Next, as a key step in the proof of Proposition 7 we will show that
\begin{equation}
\RE\ t\in [T_1,T],\ \ \sup_{m\geq 0}\sup_{n\geq 2} \E\left[\int_t^T | \tilde z^{m,n}_s|^2\ {\rm d}s\right]<+\infty.
\end{equation}
Note that for each $m\geq 0$ and $n\geq 2$, $(\tilde y_t^{m,n},\tilde z_t^{m,n})_{t\in\T}$ is an $L^1$ solution of the following BSDE:
\begin{equation}
\tilde{y}_t^{m,n}=\int_t^T \tilde{g}^{m,n}(s,\tilde{y}_s^{m,n})\ {\rm
d}s-\int_t^T \tilde{z}_s^{m,n}{\rm d}B_s,\ \ \ t\in \T,
\end{equation}
where for each $y\in \R^k$,
$\tilde{g}^{m,n}(s,y):=g^m(s,y+y_s^{m,n-1},z^{m,n-1}_s)-
g^m(s,y_s^{m,n-1},z_s^{m,n-2}).$
It follows from (H1) of $g^m$ and Remark 1 that $\as$, for each $y\in\R^k$,
\begin{equation}
\begin{array}{lll}
\left\langle y,\tilde{g}^{m,n}(t,y)\right\rangle &\leq &\left\langle y,g^m(t,y+y_t^{m,n-1},z^{m,n-1}_t)
-g^m(t,y_t^{m,n-1},z^{m,n-1}_t)
\right\rangle\vspace{0.1cm}\\
&& \Dis +|y||g^m(t,y_t^{m,n-1},z^{m,n-1}_t)
-g^m(t,y_t^{m,n-1},z^{m,n-2}_t)|\vspace{0.1cm}\\
&\leq & \Dis \kappa(|y|^2)+|y|\tilde\triangle^{m,n}_t,
\end{array}
\end{equation}
where $\kappa(\cdot)$ is defined in (H1a)$_2$, and in view of (H3) of $g^m$,
\begin{equation}
\Dis \tilde\triangle^{m,n}_t := \Dis g^m(t,y_t^{m,n-1},z^{m,n-1}_t)
-g^m(t,y_t^{m,n-1},z^{m,n-2}_t)
\leq \Dis  \lam |z^{m,n-1}_t-z^{m,n-2}_t|=\lam |\tilde z^{m,n-1}_t|
\end{equation}
as well as
\begin{equation}
\tilde\triangle^{m,n}_t\leq 2\gamma \left(g_t+
|y_t^{m,n-1}|+|z^{m,n-1}_t|+|z^{m,n-2}_t|\right)^\alpha.\vspace{0.2cm}
\end{equation}
In view of $\alpha\in (0,1/2)$ and (106), using H\"{o}lder's inequality, Jensen's inequality and Doob's inequality yields that for each $m\geq 0$ and $n\geq 2$,
\begin{equation}
\E\left[\left(\int_0^T \tilde\triangle^{m,n}_t\ {\rm d}t\right)^2\right]<+\infty\ \ {\rm and\ then}\ \
\E\left[\left.\int_0^T \tilde\triangle^{m,n}_s\ {\rm d}s\right|\F_t\right]\in \s^2.
\end{equation}
On the other hand, note from assumption (H1) of $g^m$  and (104) that $\as$,
\begin{equation}
\left\langle{\tilde{y}^{m,n}_t\over |\tilde{y}^{m,n}_t|} \mathbbm{1}_{|\tilde{y}^{m,n}_t|\neq 0},g^m(t,y^{m,n}_t,z^{m,n-1}_t)
-g^m(t,y_t^{m,n-1},z_t^{m,n-2})
\right\rangle\leq \rho(|\tilde{y}^{m,n}_t|)+\tilde\triangle_t^{m,n}.
\end{equation}
Thus, in view of (108) and (107), arguing as that from (2) to (7), by virtue of Gronwall's inequality we can obtain that for each $m\geq 0$, $n\geq 2$ and $t\in\T$,
\begin{equation}
|\tilde y_t^{m,n}|\leq \left(AT+\E\left[\left.\int_0^T \tilde\triangle^{m,n}_s\ {\rm d}s\right|\F_t\right]\right)\cdot e^{AT}
\end{equation}
and then
\begin{equation}
(\tilde y_t^{m,n})_{t\in\T}\in \s^2.
\end{equation}
In view of (110), (104), (107) and Remark 1, it follows from Proposition 1 with $p=2$ and $u=0$ that there exists a constant $C_{A,T}>0$ depending only on $A,T$ such that for each $m\geq 0$, $n\geq 2$ and $t\in \T$,
\begin{equation}
\Dis\E\left[\left(\int_t^T |\tilde z_s^{m,n}|^2\ {\rm d}s\right)\right]\leq \Dis C_{A,T}\left\{\E\left[\sup_{r\in [t,T]}|\tilde y_r^{m,n}|^2\right]+\E\left[\left(\int_t^T \tilde \triangle_s^{m,n}\ {\rm d}s\right)^2\right]+(AT)^2\right\}<+\infty.
\end{equation}
Thus, for each $m\geq 0$ and $n\geq 2$, $(\tilde y_t^{m,n},\tilde z_t^{m,n})_{t\in\T}$ is an $L^2$ solution of BSDE (103). Furthermore, in view of (104) and (105), by Proposition 3 and H\"{o}lder's inequality we can deduce that for each $m\geq 0$, $n\geq 3$ and $t\in [0,T]$,
\begin{equation}
\begin{array}{lll}
&&\Dis \E\left[\sup\limits_{r\in [t,T]}|\tilde
y^{m,n}_r|^2\right]+\E\left[\int_t^T |\tilde z^{m,n}_s|^2\ {\rm d}s\right]\vspace{0.1cm}\\
&\leq &\Dis e^{C(T-t)}\left\{\int_t^T\kappa\left(\E\left[\sup_{r\in
[s,T]}|\tilde y^{m,n}_r|^2 \right]\right){\rm d}s+\lam^2 (T-t)\E\left[\int_t^T |\tilde z^{m,n-1}_s|^2\ {\rm d}s\right]\right\},
\end{array}
\end{equation}
where the constant $C>0$ is the same as in (97).
Then, combining (112) with (98) yields that for each $m\geq 0$, $n\geq 3$ and $t\in [T_1,T]$,
\begin{equation}
\Dis \E\left[\sup\limits_{r\in [t,T]}|\tilde
y^{m,n}_r|^2\right]+\E\left[\int_t^T |\tilde z^{m,n}_s|^2\ {\rm d}s\right]
\leq \Dis 2\int_t^T\kappa\left(\E\left[\sup_{r\in
[s,T]}|\tilde y^{m,n}_r|^2 \right]\right){\rm d}s+{1\over 8}\E\left[\int_t^T |\tilde z^{m,n-1}_s|^2\ {\rm d}s\right].
\end{equation}
Note by Remark 1 that $\kappa(x)\leq A(x+1)$ for each $x\geq 0$. Gronwall's inequality together with (113) and (98) yields that for each $m\geq 0$, $n\geq 3$ and $t\in [T_1,T]$,
\begin{equation}
\begin{array}{lll}
\Dis \E\left[\sup\limits_{r\in [t,T]}|
\tilde y^{m,n}_r|^2\right]+\E\left[\int_t^T | \tilde z^{m,n}_s|^2\ {\rm
d}s\right]
&\leq &\Dis \left(2AT+{1\over 8}\E\left[\int_t^T |\tilde z^{m,n-1}_s|^2\ {\rm d}s\right]\right)\cdot e^{2A(T-t)}\vspace{0.1cm}\\
&\leq &\Dis 4AT+{1\over 4}\E\left[\int_t^T
|\tilde z^{m,n-1}_s|^2\ {\rm d}s\right]\vspace{0.1cm}\\
&\leq &\Dis {16AT\over 3}+{1\over 4^{n-2}}\E\left[\int_t^T |\tilde z^{m,2}_s|^2\ {\rm d}s\right].
\end{array}
\end{equation}
Furthermore, combining (109) and (111) with $n=2$, by virtue of Doob's inequality we obtain the existence of a constant $K_{A,T}>0$ depending only on $A$ and $T$ such that for each $m\geq 0$ and $t\in [T_1,T]$,
\begin{equation}
\E\left[\int_t^T|\tilde z^{m,2}_s|^2\ {\rm d}s\right]\leq K_{A,T}\left(1+\E\left[\left(\int_t^T \tilde \triangle_s^{m,2}\ {\rm d}s\right)^2\right]\right).
\end{equation}
And, in view of (106) with $n=2$ and the fact of  $2\alpha\in (0,1)$, it follows from H\"{o}lder's inequality and Jensen's inequality that for each $m\geq 0$ and $t\in [T_1,T]$,
$$
\begin{array}{lll}
\Dis \E\left[\left(\int_t^T \tilde \triangle_s^{m,2}\ {\rm d}s\right)^2\right]
&\leq & \Dis 16\gamma^2T^{2-2\alpha}\left(\E\left[\int_t^T g_s{\rm d}s\right]\right)^{2\alpha}+16\gamma^2T^2\E\left[\sup_{s\in [t,T]}|y^{m,1}_s|^{2\alpha}\right]\\
&&\Dis +16\gamma^2T^{2-\alpha}\E\left[\left(\int_t^T |z^{m,1}_s|^2\ {\rm d}s\right)^{\alpha}\right],
\end{array}
$$
from which as well as (67) with $n=1$, we can deduce that
\begin{equation}
\RE\ t\in [T_1,T],\ \ \sup_{m\geq 0}\E\left[\left(\int_t^T \tilde \triangle_s^{m,2}\ {\rm d}s\right)^2\right]<+\infty.
\end{equation}
Thus, the inequality (102) follows from (114)-(116).\vspace{0.1cm}

Finally, combining (101) and (102) we can deduce that for each $t\in [T_1,T]$,
\begin{equation}
\sup_{n\geq 2}\sup_{m\geq 0}\left(\E\left[\sup\limits_{r\in [t,T]}|
\bar y^{m,n}_r|^2\right]+\E\left[\int_t^T | \bar z^{m,n}_s|^2\ {\rm d}s\right]\right)<+\infty.
\end{equation}
Thus, by first taking supremum with respect to $m$ and then taking $\limsup$ with respect to $n$ in (99) as well as using Fatou's lemma, the monotonicity and continuity of the function $\kappa(\cdot)$ and Bihari's inequality, we can get (88). The proof of Proposition 7 is then completed.\vspace{0.2cm} \hfill $\Box$

Now, we come back to the proof of Theorem 4, and only consider the case when the $\alpha$ in (H3) values in $(0,1/2)$. Another case can be proved in the same way.\vspace{0.2cm}

Firstly, in view of (64) and (65) we have, for each $t\in [T_1,T]$,
\begin{equation}
\sup_{s\in [t,T]}\E\left[|y_s^m -y_s^0|\right] \leq 2\sup_{m\geq 0}\E\left[\sup_{s\in [t,T]}|y_s^{m,n}-y_s^m|\right]+\sup_{s\in [t,T]}\E\left[|y_s^{m,n} -y_s^{0,n}|\right]
\end{equation}
and for each $\beta\in (0,1)$,
\begin{equation}
\E\left[\sup_{s\in [t,T]}|y_s^m -y_s^0|^\beta\right] \leq 2\sup_{m\geq 0}\E\left[\sup_{s\in [t,T]}|y_s^{m,n}-y_s^m|^{\beta}\right]+\E\left[\sup_{s\in [t,T]}|y_s^{m,n}-y_s^{0,n}|^{\beta}\right]
\end{equation}
as well as
\begin{equation}
\Dis \E\left[\left(\int_t^T|z_s^m -z_s^0|^2{\rm d}s\right)^\beta\right]\leq \Dis 2\sup_{m\geq 0}\E\left[\left(\int_t^T|z_s^{m,n}-z_s^m|^2{\rm d}s\right)^\beta\right]
+\E\left[\left(\int_t^T|z_s^{m,n}-z_s^{0,n}|^2{\rm d}s\right)^\beta\right].
\end{equation}
Then, letting first $m\To\infty$ ($n$ being fixed) and then $n\To\infty$ in above three inequalities (118)-(120), in view of (66) and (67) in Proposition 6 and (88) in Proposition 7 we obtain that for each $t\in [T_1,T]$,
\begin{equation}
\lim\limits_{m\To \infty}\sup\limits_{s\in [t,T]}\E\left[|y_s^m-y_s^0|\right]=0,
\end{equation}
and for each $\beta\in (0,1)$,
\begin{equation}
\lim\limits_{m\To \infty}\E\left[\sup\limits_{s\in [t,T]} |y_s^m-y_s^0|^\beta+\left(\int_t^T |z_s^m -z_s^0|^2 {\rm d}s\right)^{\beta/2} \right]=0.
\end{equation}
In the sequel, set $T_j:=(T-j\delta T)\vee 0$ for each $j\geq 2$.
Noticing that the positive real number $\delta T$ depends only on $\lam$ and $A$, we can find a minimal integer $N\geq 1$ such that $T_N=0$. In view of (121) and (122), if $N=1$, then (61) and (62) have been proved. Otherwise, we consider the BSDEs with parameters $(y^m_{T_1},T_1,g^m)$ and $(y_{T_1},T_1,g)$ defined on the time interval $[T_2,T_1]$. Note by (121) that
$$\lim\limits_{m\To\infty} \E\left[|y^m_{T_1}-y_{T_1}|\right]=0.$$
Repeating the above arguments will yield that for each $t\in [T_2,T_1]$,
\begin{equation}
\lim\limits_{m\To \infty}\sup\limits_{s\in [t,T_1]}\E\left[|y_s^m-y_s^0|\right]=0,
\end{equation}
and for each $\beta\in (0,1)$,
\begin{equation}
\lim\limits_{m\To \infty}\E\left[\sup\limits_{s\in [t,T_1]} |y_s^m-y_s^0|^\beta+\left(\int_t^{T_1} |z_s^m -z_s^0|^2 {\rm d}s\right)^{\beta/2} \right]=0.
\end{equation}
Thus, in view of (121)-(124), if $N=2$, then (61) and (62) have also been proved. Otherwise, we can consider successively the BSDEs on $[T_3,T_2]$, $\cdots$, $[0,T_{N-1}]$, and finally complete the proof of Theorem 4 by repeating the above procedure.\hfill $\Box$\vspace{0.2cm}

{\bf Remark 6}\ \ The whole idea of the proof of  Theorem 4 is involved in, by virtue of (63), introducing $(y_\cdot^{m,n},z_\cdot^{m,n})$ as a bridge between $(y_\cdot^m,z_\cdot^m)$ and $(y_\cdot^0,z_\cdot^0)$, and then, by virtue of (118)-(120), proving Propositions 6 and 7 respectively. This whole idea should be  new. In addition, Propositions 6 and 7 are not easy to prove, especially a delicate argument has been done to obtain the inequality (117).\vspace{0.3cm}

In the sequel, we will investigate the stability theorem of the solutions in the space $\s^1\times {\rm M}^1$ for multidimensional BSDEs.\vspace{0.2cm}

Now, for each $m\in \N$, let $\xi^m\in \LT$ and let  $(y_t^m, z_t^m)_{t\in [0,T]}$ be
a solution in $\s^1\times {\rm M}^1$ for
the following BSDEs depending on parameter $m$:
\begin{equation}
y_t^m=\xi^m+\int_t^T g^m(s,y_s^m,z_s^m)\ {\rm d}s -\int_t^T z_s^m{\rm d}B_s,\quad t\in [0,T].
\end{equation}
Furthermore,  we introduce  the following assumptions:\vspace{0.2cm}

{\bf (B3)}\ All $g^m$ and $\xi^m$ satisfy assumptions (H1)-(H4) with the same parameters $\rho(\cdot)$, $\lambda$, $\gamma$, $g_t$ and $\alpha$. Furthermore, in the case of $\alpha\in [1/2,1)$ we assume that $g$ satisfies (H1) with a function $\rho(\cdot)$ satisfying (14) for some $\bar p>1$. \vspace{0.2cm}

{\bf (B4)}\ There exists a nonnegative real number sequence
$\{a_m\}_{m=1}^{+\infty}$ satisfying $\lim\limits_{m\To \infty}a_m=0$ such that $\as$,  for each $m\geq 1$,\vspace{-0.1cm}
\begin{equation}
\RE\ (y,z)\in \R^k\times \R^{k\times d},\ \ |g^m(\omega,t,y,z)-g^0(\omega,t,y,z)|\leq a_m.
\end{equation}
And
\begin{equation}
\lim\limits_{m\To \infty}\E\left[\sup_{t\in\T}\E\left[\left.
|\xi^{m}-\xi^{0}|\right|\F_t\right]\right]=0.\vspace{0.2cm}
\end{equation}

The following Theorem 5 is the stability theorem of the solutions in $\s^1\times {\rm M}^1$.\vspace{0.1cm}

{\bf Theorem 5}\ Under assumptions (B3) and (B4), we have
\begin{equation}
\lim\limits_{m\To \infty}\E\left[\sup\limits_{t\in [0,T]} |y_t^m-y_t^0|+\left(\int_0^T |z_s^m -z_s^0|^2\  {\rm d}s\right)^{1/2} \right]=0.\vspace{0.1cm}
\end{equation}

{\bf Proof}\ \ First of all, it follows from (62) of Theorem 4 that the sequence of variables $\{\sup_{t\in [0,T]} |y_t^m-y_t^0|\}_{m=1}^{+\infty}$ converges in probability to $0$ as $m\To\infty$. In the sequel, we will first prove that it is also uniformly integrable, and then
\begin{equation}
\lim\limits_{m\To \infty}\E\left[\sup\limits_{t\in [0,T]} |y_t^m-y_t^0|\right]=0.\vspace{0.2cm}
\end{equation}
In fact, set $\hat y^m_\cdot:=y^m_\cdot-y^0_\cdot$ and $\hat z^m_\cdot:=z^m_\cdot-z^0_\cdot$. Then, for each $m\geq 1$, $(\hat y^m_t,\hat z^m_t)_{t\in \T}$ is a solution in the space $\s^1\times {\rm M}^1$ for the following BSDE:
\begin{equation}
\hat{y}_t^m=\xi^{m}-\xi^0+\int^{T}_{t}\hat g^m(s,\hat y_s^m,\hat z_s^m)\ \mathrm{d} s-\int^{T}_{t}\hat{z}_s^m\ \mathrm{d}B_s, \quad t\in\T,
\end{equation}
where for each $y\in\R^k$,
$\hat g^m(t,y,z):=g^m(t,y+y^0_t,z+z^0_t)
-g^0(t,y^0_t,z^0_t).$
It follows from assumptions (H1) and (H3) of $g^m$ together with (126) and Remark 1 that $\as$, for each $y\in\R^k$,
\begin{equation}
\begin{array}{lll}
\langle y,\hat{g}^m(t,y,z)\rangle &\leq &\langle y,g^m(t,y+y_t^0,z+z^0_t)-g^m(t,y_t^0,z+z^0_t)
\rangle+|y||g^m(t,y_t^0,z+z^0_t)-g^m(t,y_t^0,z^0_t)|\\
&& \Dis +|y||g^m(t,y_t^0,z^0_t)-g^0(t,y_t^0,z^0_t)|\\
&\leq & \Dis \kappa(|y|^2)+\lam |y||z|+|y|a_m
\end{array}
\end{equation}
and
\begin{equation}
\left\langle{\hat{y}^m_t\over |\hat{y}^m_t|} \mathbbm{1}_{|\hat{y}^m_t|\neq 0},g^m(t,y^m_t,z^m_t)
-g^0(s,y_t^0,z_t^0)\right\rangle\leq \rho(|\hat{y}^m_t|)+\hat\triangle_t^m+a_m,\vspace{0.2cm}
\end{equation}
where $\kappa(\cdot)$ is defined in (H1a)$_2$, and
\begin{equation}
\hat\triangle^m_t:=|g^m(t,y_t^0,z^m_t)-g^m(t,y_t^0,z^0_t)| \leq 2\gamma \left(g_t+
|y_t^0|+|z^m_t|+|z^0_t|\right)^\alpha.
\end{equation}
Next, in view of (132) and Remark 1, using Corollary 2.3 in \cite{Bri03} and Gronwall's inequality we get that for each $m\geq 1$ and $t\in\T$,
\begin{equation}
|\hat y_t^m|\leq \left(\E\left[\left.|\xi^m-\xi^0|\right|\F_t
\right]+AT+\E\left[\left.\int_0^T \hat\triangle^m_s\ {\rm d}s\right|\F_t\right]+a_mT\right)\cdot e^{AT}.
\end{equation}
On one hand, it follows from (127) that $\{\sup_{t\in\T}\E\left[\left.
|\xi^{m}-\xi^{0}|\right|\F_t\right]\}_{m=1}^{+\infty}$ is uniformly integrable. On the other hand, in view of (133), for each $q>1$ such that $\alpha q<1$, H\"{o}lder's inequality yields the existence of a positive constant $C$ independent of $m$ such that
\begin{equation}
\E\left[\left(\int_0^T \hat\triangle^m_s\ {\rm d}s\right)^q\right]<C\left(1+\E\left[\left(\int_0^T |z^m_s|^2\ {\rm d}s\right)^{\alpha q\over 2} \right]\right).
\end{equation}
And, Doob's inequality together with (135) and (62) leads to
\begin{equation}
\sup_{m\geq 1}\E\left[\sup_{t\in\T}\left( \E\left[\left.\int_0^T \hat\triangle^m_s\ {\rm d}s\right|\F_t\right]\right)^q\right]\leq
\sup_{m\geq 1}\E\left[\left(\int_0^T \hat\triangle^m_s\ {\rm d}s\right)^q\right]<+\infty,
\end{equation}
which means that the sequence of variables $$\left\{\sup_{t\in\T}\E\left[\left.\int_0^T \hat\triangle^m_s\ {\rm d}s\right|\F_t\right]\right\}_{m=1}^{+\infty}$$
is uniformly integrable. Thus, in view of (134) and (B4), we can deduce that the sequence of variables $\{\sup_{t\in [0,T]} |\hat y_t^m|\}_{m=1}^{+\infty}$ is also uniformly integrable, and get (129).\vspace{0.1cm}

Finally, in view of (131), Remark 1, Lemma 1 and (129), using Proposition 1 with $p=1$ and a similar argument to that from (30) to (33) we can obtain that \begin{equation}
\lim\limits_{m\To \infty}\E\left[\left(\int_0^T |z_s^m -z_s^0|^2\  {\rm d}s\right)^{1/2} \right]=0.
\end{equation}
Thus, Theorem 5 is proved by (129) and (137).\vspace{0.3cm} \hfill $\Box$

{\bf Remark 7}\ We mention that Theorem 5 can also be proved in the same way as Theorem 4. In fact, it can be proved that (67) holds for $\beta=1$ under the assumptions (B1) and (B2) of Theorem 5.\vspace{0.2cm}

By Theorems 4 and 5, the following corollary follows immediately.\vspace{0.1cm}

{\bf Corollary 2}\ \ Assume that the generator $g$ satisfies assumptions (H1)-(H3). In the case when the $\alpha$ defined in (H3) values in $[1/2,1)$, we also assume that there exists a constant $\bar p>1$ such that the function $\rho(\cdot)$ in (H1) satisfies (14). Suppose that for each $m\geq 1$, $\xi^m,\xi\in \LT$ and $(y^m_t,z^m_t)_{t\in\T}$ and $(y_t,z_t)_{t\in\T}$ are respectively the unique $L^1$ solution of BSDE$(\xi^m,T,g)$ and BSDE$(\xi,T,g)$. If $\lim\limits_{m\To\infty}\E\left[|\xi^m-\xi|\right]=0,$ then
$$
\lim\limits_{m\To \infty}\sup\limits_{t\in [0,T]}\E\left[|y_t^m-y_t|\right]=0,
$$
and for each $\beta\in (0,1)$,
$$
\lim\limits_{m\To \infty}\E\left[\sup\limits_{t\in [0,T]} |y_t^m-y_t|^\beta+\left(\int_0^T |z_s^m -z_s|^2 {\rm d}s\right)^{\beta/2} \right]=0.
$$
Moreover, if for each $m\geq 1$, $\xi^m$ and $g(t,0,0)$ satisfy assumptions (H4), and
$$
\lim\limits_{m\To \infty}\E\left[\sup_{t\in\T}\E\left[\left.
|\xi^{m}-\xi|\right|\F_t\right]\right]=0,
$$
then
$$
\lim\limits_{m\To \infty}\E\left[\sup\limits_{t\in [0,T]} |y_t^m-y_t|+\left(\int_0^T |z_s^m -z_s|^2\  {\rm d}s\right)^{1/2} \right]=0.\vspace{0.3cm}
$$

{\bf Remark 8}\ By Remarks 5 and 2, we know that Theorems 4-5 and Corollary 2 give the stability of the $L^1$ solutions of multidimensional BSDEs investigated in \citet{Bri03}.\vspace{0.1cm}




\end{document}